\documentclass{gtmon_a}     % Basic GT/GTM/AGT style
\pdfoutput=1

\usepackage{graphicx}  %%% the recommended graphics package
\usepackage{pinlabel}  %%% the recommended graphics+labelling package

\proceedingstitle{Heegaard splittings of 3-manifolds (Haifa 2005)}
\conferencestart{10 July 2005}
\conferenceend{19 July 2005}
\conferencename{Heegaard splittings of 3-manifolds}
\conferencelocation{Haifa}

\editor{Cameron Gordon}
\givenname{Cameron}
\surname{Gordon}

\editor{Yoav}
\givenname{Yoav}
\surname{Moriah}

\title{Incompressible surfaces and $(1,2)$--knots}

\author[M\,Eudave-Mu\~noz]{Mario Eudave-Mu\~noz}
\givenname{Mario}
\surname{Eudave-Mu\~noz}
\address{Instituto de Matem\'aticas\\ Universidad Nacional Aut\'onoma de M\'exico\\\newline
Ciudad Universitaria\\ 04510 M\'exico D.F.\\ MEXICO\\\newline 
and\\Centro de Investigaci\'on en
Matem\'aticas\\\newline  Apdo. Postal 402\\ 36000 Guanajuato, Gto.\\ MEXICO}
\email{mario@matem.unam.mx}
\urladdr{}


\arxivreference{math.GT/0703132}

\volumenumber{12}
\issuenumber{}
\publicationyear{2007}
\papernumber{03}
\startpage{35}
\endpage{87}

\doi{}
\MR{}
\Zbl{}

\received{4 October 2006}
\revised{7 February 2007}
\accepted{6 April 2007}
\published{3 December 2007}
\publishedonline{3 December 2007}
\proposed{}
\seconded{}
\corresponding{}
\version{}

\makeatletter
\def\cnewtheorem#1[#2]#3{\newtheorem{#1}{#3}[section]
\expandafter\let\csname c@#1\endcsname\c@thm}

\AtBeginDocument{\let\bar\wbar\let\tilde\wtilde\let\hat\what}
\def\S{Section }
\def\Int{\mathop{\mathrm{int}}\nolimits}
\makeop{genus}
\makeop{tn}

\newtheorem{thm}{Theorem}[section]    % Standard theorem environment
\cnewtheorem{cla}[thm]{Claim}   
\makeautorefname{cla}{Claim}
\makeatother  %  move after \newtheorem block

\begin{document}

\begin{asciiabstract}
We give a description of all (1,2)-knots in S^3
which admit a closed meridionally incompressible surface of genus 2 in their complement.
That is, we give several constructions of (1,2)-knots having a meridionally incompressible
surface of genus 2, and then show that any such surface for a 
(1,2)-knot must come
from one of the constructions. As an application, we show explicit examples of tunnel number
one knots which are not (1,2)-knots.
\end{asciiabstract}

\begin{htmlabstract}
We give a description of all (1,2)&ndash;knots in S<sup>3</sup>
which admit a closed meridionally incompressible surface of genus 2 in their complement.
That is, we give several constructions of (1,2)&ndash;knots having a meridionally incompressible
surface of genus 2, and then show that any such surface for a
(1,2)&ndash;knot must come
from one of the constructions. As an application, we show explicit examples of tunnel number
one knots which are not (1,2)&ndash;knots.
\end{htmlabstract}

\begin{abstract}
We give a description of all $(1,2)$--knots in $S^3$
which admit a closed meridionally incompressible surface of genus 2 in their complement.
That is, we give several constructions of $(1,2)$--knots having a meridionally incompressible
surface of genus 2, and then show that any such surface for a 
$(1,2)$--knot must come
from one of the constructions. As an application, we show explicit examples of tunnel number
one knots which are not $(1,2)$--knots.
\end{abstract}

\maketitle

%%%%%%%%%%%%%%%%%%%%   Start of main body of article

\section{Introduction}

An important problem in knot theory is that of determining all 
incompressible surfaces in a given knot complement.  
The main purpose of this paper is to give a description of all $(1,2)$--knots in $S^3$
which admit a closed meridionally incompressible surface of genus 2 in their complement.
Another purpose is to construct explicit examples of tunnel number one knots 
that are not  $(1,2)$--knots. These are obtained by combining the constructions
of  the author \cite{E2} with the results of this paper.

Let $F$ be a closed surface of genus $g$ standardly embedded in $S^3$,
that is, it bounds a handlebody on each of its sides. Following Doll \cite{D}, we say
that a knot $K$ has a $(g,b)$--{\it presentation\/} or that it is in a $(g,b)$--{\it position\/}, if 
$K$ has been isotoped to intersect $F$
transversely in $2b$ points that divide $K$ into $2b$ arcs, so that the
$b$ arcs in each side can be isotoped, keeping the endpoints fixed, to
disjoint arcs on $F$. 
The {\it genus--g--bridge number\/} of $K$, $b_g(K)$, is the smallest integer $n$ for which $K$ has a 
$(g,n)$--presentation.
The genus--0--bridge number $b_0(K)$ is then the usual bridge number; we say that a knot is an $n$--bridge knot if $b_0(K)\leq n$. 
Here we will consider only the case $g=1$. We say that a knot is a $(1,n)$--{\it knot\/} if $b_1(K)\leq n$.
It is not difficult to see that if $K$ has a $(g,b)$--presentation, then the tunnel number of $K$, denoted
$\tn(K)$, satisfies $\tn (K)\leq g+b-1$.

Let $K$ be a knot in $S^3$, and $S$ a closed surface in its complement. 
We say that $S$ is meridionally compressible  
if there is an embedded disk $D$ in $S^3$, with $S\cap D=\partial D$ 
 a nontrivial curve in $S$, and so that $K$ intersects $D$ 
at most in one point. Otherwise $S$ is called meridionally incompressible.
In particular, if $S$ is meridionally incompressible then it is
incompressible in $S^3-K$. 

Incompressible surfaces in the complement of knots with a given bridge number have been studied in 
several cases. Schubert \cite{Sc} studied incompressible tori in the complement of knots and found a 
relation between the bridge numbers of a satellite knot $K$ and its companion $K'$. 
Hatcher and Thurston \cite{HT} proved that 
there are no closed incompressible surfaces in the complement of $2$--bridge knots; this also follows
from Gordon and Litherland \cite{GLi}. On the other hand, \mbox{Finkelstein and Moriah \cite{FM}} and \mbox{Wu \cite{W}} proved that 
the complement of a {\it generic\/} $n$--bridge knot, $n\geq 3$, contains a closed incompressible surface 
(but in general the surface is meridionally compressible). More recently Ozawa \cite{O} has given a
description of all $3$--bridge knots whose complement contain a closed incompressible surface of genus 2.
It may be difficult to do something similar for $3$--bridge knots and surfaces of higher genus.
In Eudave-Mu{\~n}oz and Neumann-Coto \cite{EN}, some examples are given of $3$--bridge knots whose complement contain meridionally
incompressible surface of arbitrarily high genus, ie examples of a single $3$--bridge knot which 
contains infinitely many closed meridionally incompressible surfaces in its complement.

Incompressible surfaces in the complement of $(1,1)$--knots have also been studied. 
$(1,1)$--knots whose complement contain an incompressible torus, ie satellite $(1,1)$--knots, have been
classified by Morimoto and Sakuma \cite{MS}; see also Eudave-Mu{\~n}oz~\cite{E1}. Well, in those papers satellite tunnel
number one knots are classified but these turn out to be $(1,1)$--knots, as it is shown in \cite{E1}. 
In \cite{E2} a  construction is given of $(1,1)$--knots whose complement contain a closed meridionally
incompressible surface of genus $g$, and in \cite{E4} it is proved that any $(1,1)$--knot whose
complement contains a closed meridionally incompressible surface must come from that construction. 
It is shown by Saito \cite{Sa} that the complement of a satellite tunnel number one knot does not 
contain any  meridionally incompressible surface other than the satellite torus; this implies that the
knots constructed in \cite{E2} are hyperbolic.
It follows from work of Gordon and Reid \cite{GR} that the complement of a
$(1,1)$--knot cannot contain an incompressible planar meridional surface, ie a meridional surface 
of genus 0 (a meridional surface is a properly embedded surface in a knot exterior whose boundary
consists of meridians of the knot). On the other hand, in Eudave-Mu{\~n}oz and Ram\'{\i}rez-Losada \cite{ER1} a description is given of all
$(1,1)$--knots whose complement contain a meridional and meridionally incompressible surface of 
genus $g\geq 1$. The complement of any of these knots also contains a closed incompressible surface (but
perhaps meridionally  compressible) by Culler, Gordon, Luecke and Shalen \cite{CGLS}.

Not much is known about incompressible surfaces in the complement of $(1,2)$--knots. In \cite{EN}, a
construction is given of hyperbolic $(1,2)$--knots whose complement contain an acylindrical surface
of genus $g$, $g\geq 2$, ie  an incompressible
surface which divides the exterior of the knot into manifolds that do not contain essential annuli.
The example given in \cite{EN} of a $3$--bridge knot whose complement contains meridionally incompressible
surfaces of arbitrarily high genus can be adapted to produce an example of a $(1,2)$--knot whose
complement contains meridionally incompressible surfaces of arbitrarily high genus. 
To do that just embed the branched surface given in \cite[Figure~14]{EN} as a surface of type 6, 
defined in \fullref{sec:type6} of this paper. This example shows that it may be difficult to give a 
description of all $(1,2)$--knots whose complement contain a closed meridionally incompressible
surface. However, in this paper we do this for surfaces of genus $2$. In \fullref{sec:construction} 
we give several constructions which produce $(1,2)$--knots whose complement contain
a closed meridionally incompressible surface of genus $2$. In 
\fullref{sec:characterization} we show that if the 
complement of a $(1,2)$--knot contains a closed meridionally incompressible surface of genus 2,
then the knot and the surface come from one of the given constructions. Contained in that proof is
also a description of all $(1,2)$--knots whose complement contain a meridionally incompressible
torus; these are some satellites of $(1,1)$--knots. It also follows from the construction that there are
$(1,2)$--knots whose complement contains both a closed meridionally incompressible surface of
genus 2 and of genus 1.

If $K$ is a $(1,1)$--knot, then it is easy to see that $K$ has tunnel number one. On the other hand,
if $K$ has tunnel number one, it seems to be very difficult to determine $b_1(K)$. A priori 
there should be tunnel number one knots with arbitrarily large $b_1$, but this has been 
difficult to prove. Moriah and Rubinstein \cite{MR} showed the existence of tunnel number one knots $K$ with
$b_1(K)\geq 2$. Morimoto, Sakuma and Yokota also showed this, and gave explicit examples of 
knots $K$ with tunnel number one and $b_1(K)=2$ \cite{MSY}. It was shown in \cite{E4} 
that many of the tunnel  number one knots $K$ constructed in \cite{E2}  are not $(1,1)$--knots. We refine
the proof here, and show in \fullref{sec:not1-2knots} that some of the knots constructed in \cite{E2} are not
$(1,2)$--knots. The argument is as follows: the complement of the tunnel number one knots constructed 
in \cite{E2} contains a closed
meridionally incompressible surface. We then pick some of them whose complement contain a meridionally
incompressible surface of genus 2, and show that the surface comes from none of the constructions of
\fullref{sec:construction}. This implies that $b_1(K)\geq 3$ for any such knot $K$. These knots can be explicitly
constructed;  an example is given in \fullref{fig:not1-2knot}. However, our examples seem to satisfy 
 $b_1(K)\geq 4$,  and we do not know if one of these knots satisfies $b_1(K)=3$.

Recently, Johnson and Thompson \cite{JT}, and independently Minsky, Moriah and Schlei\-mer \cite{MMS} have
shown that for any given $n$, there exist tunnel number one knots which are not $(1,n)$--knots.
In \cite{MMS} it is in fact shown that for given $t$ and $n$, there exist tunnel number $t$ knots $K$
so that $b_t(K) > n$. The two papers use similar techniques and prove the existence of such knots, 
but do not give explicit examples.
 
Finally, we add that Valdez-S\'anchez and Ram\'{\i}rez-Losada \cite{RV2} have also shown examples of 
tunnel number one knots $K$ with $b_1(K)=2$. These knots bound punctured Klein bottles but are not
contained in the $(1,1)$--knots bounding Klein bottles determined by the same authors \cite{RV1}.
In Eudave-Mu\~noz \cite{E3} a construction is given of tunnel number one knots which admit a meridional 
incompressible surface. Using \cite{ER1}, we can show that among these knots there are ones 
with $b_1(K)=2$ \cite{ER2}, and expect to prove that there are others with $b_1(k)=3$.

\section{Construction of meridionally incompressible surfaces}\label{sec:construction}

Let $F$ be a closed surface of genus $g$ standardly embedded in $S^3$,
and let $I=[0,1]$.
Consider a product neighborhood $F\times I$ of $F$. To say that a knot $K$ has a 
$(g,b)$--presentation is equivalent to say that $K$ has been isotoped to lie in $F\times I$, so that
$K\cap (F\times \{ 0 \})$ and $K\cap (F\times \{ 1 \})$ consists each of
$b$ arcs (or $b$ tangent points), and the rest of the knot consists of $2b$ vertical arcs 
in $F\times I$, that is, arcs which intersect each leave $F\times \{t\}$, $0< t< 1$, in the
product exactly in one point. Or simply, $K$ is in a $(g,b)$--position if $K \subset F\times I$, 
and the projection map
$p\co F\times I \rightarrow I$ when restricted to $K$ has exactly $b$ local maxima and $b$ local minima.

Let $T$ be a standard torus in $S^3$, and let $I=[0,1]$.  Consider
$T\times I\subset S^3$.
$T_0=T\times \{ 0\}$ bounds a solid torus $R_0$, and $T_1=T\times \{1\}$ bounds a
solid torus $R_1$, such that $S^3= R_0\cup (T\times I)\cup R_1$.
Think of the solid torus $R_1$ as containing the point at infinity.
By a vertical arc in a product $T\times [a,b]$ we mean an
embedded arc which intersects every torus $T\times \{x\}$  in
the product in at most one point. By a level simple closed curve we mean a curve 
which is contained in some level torus $T\times \{ x\}$.

In this section we construct knots with a $(1,2)$--presentation whose complement contain a closed
meridionally incompressible surface of genus $2$.
We assume that all knots constructed in this section are contained
in $T\times I$.

\subsection{Surfaces of type 1}\label{sec:type1}

Choose a point $e$ on $I$, so that
$0 < e < 1$. Consider the torus $T_e = T\times \{e\}$.
Let $\gamma_0$, $\gamma_1$ be simple closed curves embedded in the
product $T\times (0,e)$ and $T\times (e,1)$ respectively, so that each
curve has only one local maximum and one local minimum with respect to the
projection to $(0,e)$ or to $(e,1)$. Suppose also that  $\gamma_0$ ($\gamma_1$) 
is not in  a $3$--ball contained in
$R_0\cup (T\times [0,e])$ (resp.\ $(T\times [e,1])\cup R_1$), that is, 
it is not a trivial knot in that region.

Let $\alpha$ be a vertical arc in
$T\times [0,1]$, joining the maximum point of $\gamma_0$ with the
minimum of  $\gamma_1$. Let $\Gamma_1$ be the $1$--complex consisting of
the union of  the curves $\gamma_0$, $\gamma_1$ and the arc $\alpha$. So
$\Gamma_1$ is a trivalent graph embedded in $S^3$. 

\begin{figure}[ht!]
\labellist
\small
\hair 2pt
\pinlabel $T_1$ at 7 123
\pinlabel $T_0$ at 69 69
\pinlabel $\gamma_1$ at 19 61
\pinlabel $\gamma_0$ at 90 93
\pinlabel $\alpha$ at 48 95 
\pinlabel $\gamma_1$ at 315 90
\pinlabel $\gamma_0$ at 217 113
\pinlabel $\alpha$ at 257 52
\endlabellist
\centering
\includegraphics{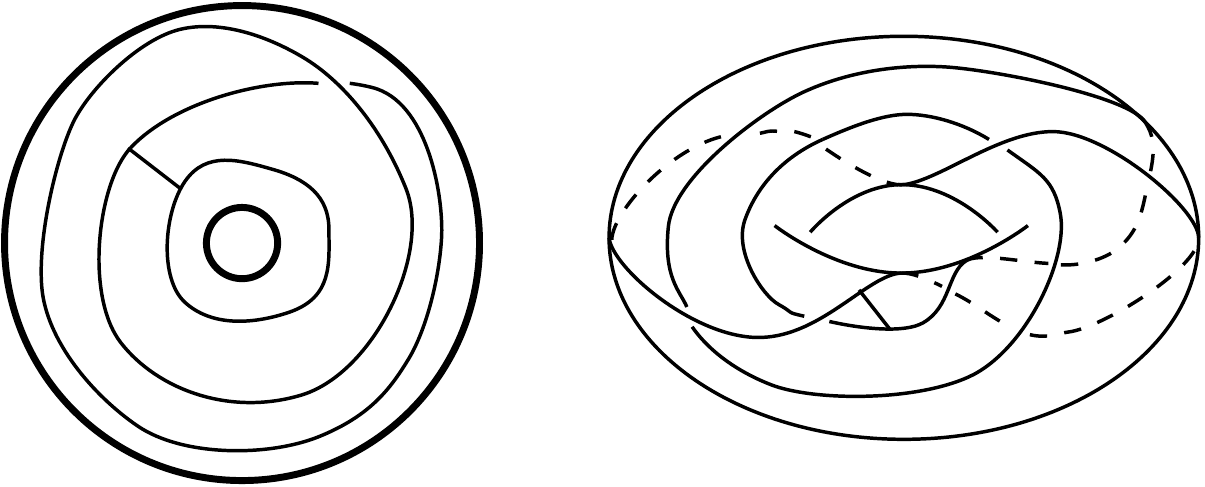}
\caption{}
\label{fig:type1}
\end{figure}

Let $N(\Gamma_1)$ be a regular neighborhood of
$\Gamma_1$. This is a genus $2$ handlebody. We can assume that 
$N(\Gamma_1)$
is the union of $2$ solid tori, $N(\gamma_0)$ and $N(\gamma_1)$, joined
by the $1$--handle $N(\alpha)$. Let $S=\partial N(\Gamma_1)$, we say that $S$ is a surface of type 1.
In \fullref{fig:type1} we show in a schematic way a surface of type 1, and also give an explicit example,
where both curves $\gamma_i$, $i=1,0$, are isotopic to level curves.

There are three possibilities for the graph $\Gamma_1$:

(1)\qua By isotopies of $\Gamma_1$, none of the curves $\gamma_0$ and $\gamma_1$ 
can be isotoped into a level curve.

(2)\qua By an isotopy of $\Gamma_1$, both curves $\gamma_0$ and $\gamma_1$ can be isotoped 
into level curves, say $\gamma_i$ can be isotoped into $\gamma_i'$, which lies in $T_i$, $i=0,1$. 
In this case we assume that $\gamma_i'$ is not isotopic to the core of $R_i$, ie
$\gamma_i'$ does not consist of a longitude and several meridians of $R_i$, and assume also that
$\Delta(\gamma_0',\gamma_1')\geq 2$ \cite{E2}.

(3)\qua By an isotopy of $\Gamma_1$, only one of the curves, say $\gamma_1$, can be isotoped into
a level curve $\gamma_1'$ contained in $T_1$. Assume that $\gamma_1'$ is not isotopic to the core of
$R_1$, ie $\gamma_1'$ does not consist of a longitude and several meridians of $R_1$.

\begin{figure}[ht!]
\centering
\includegraphics{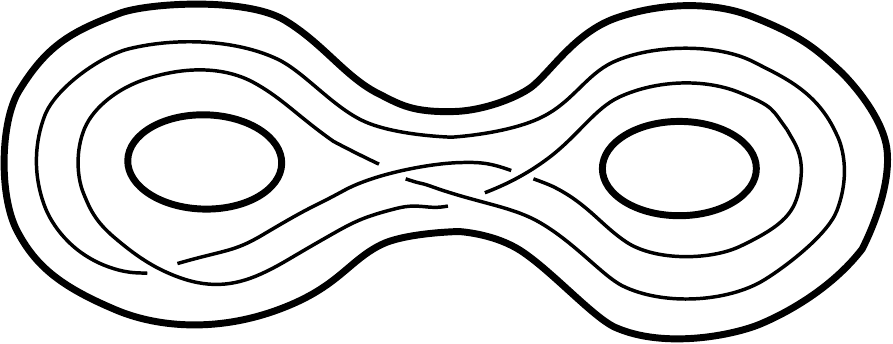}
\caption{}\label{fig:type1s}
\end{figure}

Let $E=N(\Gamma_1)\cap T_e$, this is a disk, which is in fact a cocore of the 
$1$--handle $N(\alpha)$. Embed a knot $k$ in $N(\Gamma_1)$ so that it intersects $E$ in 
four points, $k\cap N(\gamma_1)$ consists of two arcs each having just one local maximum,
$k\cap N(\gamma_0)$ consists of two arcs each having just one local minimum, and $k\cap N(\alpha)$
consists of four vertical arcs.
Suppose also that $S=\partial N(\Gamma_1)$ is meridionally incompressible in $N(\Gamma_1)-k$. 
It is not difficult 
to see that there are plenty of such knots. See \fullref{fig:type1s} for an example.

\begin{thm}\label{thm:type1} Let $\Gamma_1$ be a graph as above, and $k$ in $N(\Gamma_1)$ as above. 
Then $k$ is a $(1,2)$--knot, and $S=\partial N(\Gamma_1)$ is  meridionally incompressible in $S^3-k$,
except possibly if $\Gamma_1$ is as in case (3), ie it can be isotoped so that exactly one
of the curves $\gamma_i$ is a level curve.
\end{thm}

\begin{proof} By construction $k$ has a $(1,2)$--presentation, and by hypothesis $S$ 
is meridionally incompressible in $N(\Gamma_1)-k$. So it remains to prove that $S$ is
incompressible in $S^3-\Int N(\Gamma_1)$.

Let $T_e' = T_e - \Int N(\Gamma_1)$, this is a once-punctured torus.
Suppose $D$ is a compression disk for $S$, and suppose it intersects
transversely the torus $T_e'$. Let $\beta$ be a simple closed curve of
intersection between $D$ and $T_e'$, which is innermost
in $D$. So $\beta$ bounds a disk $D' \subset D$, which is contained, say, in the
solid torus $(T\times [e,1])\cup R_1$. If $\beta$
is trivial on $T_e'$, then by cutting $D$ with an innermost disk lying in the disk 
bounded by $\beta$ on $T_e'$, we get a compression disk with fewer intersections with 
$T_e'$. If $\beta$ is essential on $T_e'$, then it would be parallel to
$\partial T_e'$, or it would be a meridian of the solid torus $(T\times [e,1])\cup R_1$, 
but in any case the curve $\gamma_1$ will be contained in a $3$--ball, which is a
contradiction.

So suppose $D$ intersects $T_e'$ only in arcs. Let $\beta$ be such
an arc which is outermost on $D$; it cobounds with an arc
$\delta \subset \partial D$ a disk $D'$. If $\beta$ is parallel to an arc
on $\partial T_e'$, then by cutting $D$ with such an outermost arc lying
on $T_e'$ we get another compression disk with fewer intersections with
$T_e'$, so assume that $\beta$ is an essential arc on $T_e'$.
After isotoping $D$ if necessary, we can assume that the arc $\delta$ can
be decomposed as
$\delta = \delta_1 \cup \delta_2 \cup \delta_3$, where
$\delta_1, \delta_3$ lie on $\partial N(\alpha)$ and $\delta_2$
lie on $\partial N(\gamma_1)$ (if $\delta$ were contained in
$\partial N(\alpha)$, then by isotoping $D$ we would get a compression
disk whose intersection with $T_e'$ contains a simple closed curve). Let
$E$ be a disk contained in $\partial N(\alpha)$ so that
$\partial E = \delta_1 \cup \delta_4 \cup \delta_3 \cup \delta_5$, where
$\delta_4$ lies on $\partial T_e'$ and $\delta_5$ lies on 
$\partial N(\gamma_1)\cap \partial N(\alpha)$.
So $D' \cup E$ is an annulus, where one boundary component, ie
$\beta \cup \delta_4$ lies on $T_e$, and the other,
$\delta_2 \cup \delta_5$, lies on $\partial N(\gamma_1)$. If $\delta_2
\cup \delta_5$ is a meridian of $\gamma_1$, then $\beta \cup \delta_4$ 
is a meridian of the solid torus $(T\times [e,1])\cup R_1$. Then
$\gamma_1$ intersects a meridian disk of $(T\times [e,1])\cup R_1$ in one
point, which implies that it is parallel to a knot lying on the torus
$T_1$, and it is isotopic to the core of $R_1$, which is a contradiction. 
If $\delta_2 \cup \delta_5$ goes
more than once longitudinally on $\gamma_1$, then the level curve
$\beta\cup \delta_4$ (which is a trivial or a torus knot) would be a 
cable knot around $\gamma_1$. This shows that $\gamma_1$
is a core of the solid torus $R_1$, which is not
possible.  If $\delta_2 \cup \delta_5$ goes
once longitudinally on $\gamma_1$, then $\gamma_1$ is
isotopic to a curve on $T_1$.
So we conclude that either $S$ is incompressible or $\gamma_1$ is
isotopic to a level curve on $T_1$.

Suppose now that $\Gamma_1$ has been isotoped so that $\gamma_1$ is a level curve in $T_1$
and that $\gamma_0$ is a level curve in $T_0$.
The incompressibility of $S$ now follows from \cite[Theorem 4.1]{E2}, because
we assume that $\Delta(\gamma_1,\gamma_2) \geq 2$. This completes the proof.
\end{proof}

Note that if one of the curves $\gamma_0$, $\gamma_1$ does not satisfy the required conditions,
then the surface $S$ will be compressible.

It is not difficult to construct examples where $\gamma_1$ is a level curve, but 
the curve $\gamma_0$ it is not, and so that $S$ is incompressible in $S^3-\Int N(\Gamma_1)$.
But at this writing we do not have a precise description of all such curves.
On the other hand, it is also not difficult to construct examples of graphs with such curves so 
that the surface $S$ is in fact compressible. One such example can be constructed starting with 
a graph $\Gamma_1$ so that $S=\partial N(\Gamma_1)$ is obviously compressible, say a graph where
both curves $\gamma_0$ and $\gamma_1$ are level and $\Delta(\gamma_0,\gamma_1)=1$. Now slide an 
endpoint of $\gamma_0$ through $\alpha$ and then through $\gamma_1$, going around it several times, 
and then again through $\alpha$, to get a new curve $\gamma_0'$ and a new graph $\Gamma_1'$. 
The new curve $\gamma_0'$ can be chosen so that it has a single local maximum, a local minimum 
and it is not isotopic to a level curve. However $S'=\partial N(\Gamma_1')$ would be
compressible, because the exteriors of $\Gamma_1$ and $\Gamma_1'$ are homeomorphic.

\subsection{Surfaces of type 2}\label{sec:type2}

Let $e$ be a point on $I$, so that
$0 < e < 1$. Consider the level torus $T_e = T\times \{e\}$.
Let $\gamma$ be a simple closed curve embedded in the
level torus $T_e$. Let $\alpha$ be an arc contained in $T\times I$,
with endpoints in $\gamma$, so that it has just a local maximum at $T_1$, and 
just a local minimum at $T_0$. Suppose that $\gamma$ is essential in $T_e$, or well, it is 
inessential but bounds a disk in $T_e$ which intersects $\alpha$ in one point.
Note that the interior of $\alpha$ intersects $T_e$ in one point, so $\alpha$ 
is divided into a lower and an upper arc, say $\alpha_1$ and $\alpha_2$. 
Suppose that none of these arcs can be isotoped (in $R_0\cup (T\times [0,e])$
or $(T \times [e,1])\cup R_1$),
keeping its endpoints fixed, into an arc on $T_e$ with interior disjoint from $\gamma$.

\begin{figure}[ht!]
\labellist
\small\hair 2pt
\pinlabel $T_1$ at 9 123
\pinlabel $T_0$ at 69 69
\pinlabel $\gamma$ at 29 86
\pinlabel $\alpha$ at 83 88
\pinlabel $\gamma$ at 265 84
\pinlabel $\alpha$ at 244 18
\endlabellist
\centering
\includegraphics{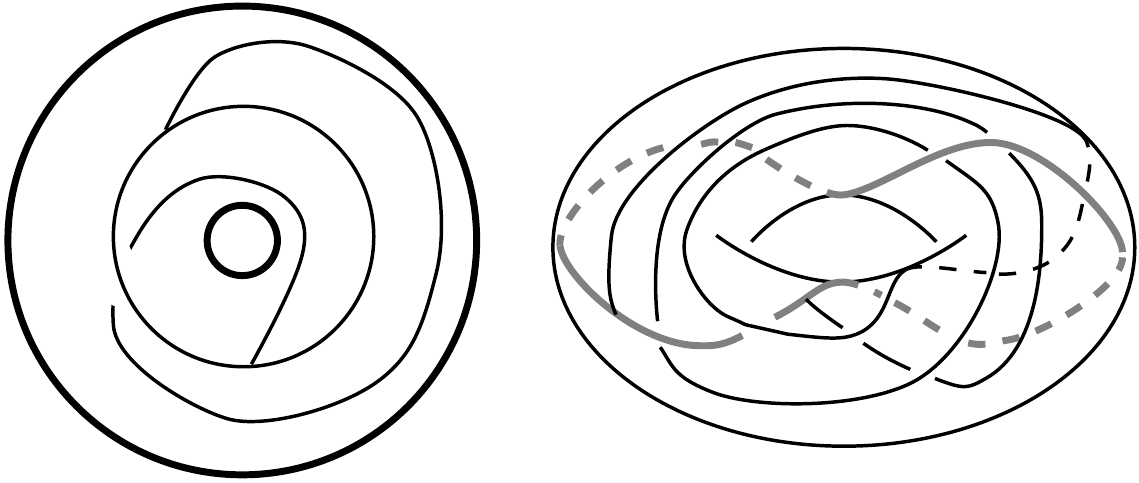}
\caption{}\label{fig:type2}
\end{figure}

Let $\Gamma_2$ be the $1$--complex consisting of
the union of  the curve $\gamma$ and the arc $\alpha$. So
$\Gamma_2$ is a trivalent graph embedded in $S^3$. 
Let $N(\Gamma_2)$ be a regular neighborhood of
$\Gamma_2$. This is a genus $2$ handlebody. We can assume 
that $N(\Gamma_2)$
is the union of the solid tori $N(\gamma)$ and $1$--handle $N(\alpha)$, so that these intersect in two disks,
$E_1$  and $E_2$, where say $E_1$ is at level $T\times \{ e-\epsilon\}$ and $E_2$ at level 
$T\times \{e+\epsilon\}$, for some small $\epsilon > 0$.
Let $S=\partial N(\Gamma_2)$. We say that $S$ is a surface of type 2.
In \fullref{fig:type2} we show schematically a surface of type 2, and give an explicit example.

Let $k$ be a knot embedded in $N(\Gamma_2)$ so that $k$ intersects each of 
$E_1$ and $E_2$ in two points, $k\cap N(\alpha)$ consists of two arcs, each
with a local maximum and a local minimum in $T\times I$, and $k\cap N(\gamma)$,
which is contained in $T\times [e-\epsilon,e+\epsilon]$, consists of two 
vertical arcs. Suppose also that $S=\partial N(\Gamma_2)$ is meridionally incompressible in
$N(\Gamma_2)-k$.  To get that it suffices to ask that $k$ is well wrapped in $N(\Gamma_2)$ (ie consider
the two arcs of $k$ lying in $N(\gamma)$, get a knot by joining the ends of the arcs 
lying in $E_1$ and $E_2$ with an arc contained in such disks, and then push the knot to the interior of
$N(\gamma)$;  to be well wrapped just means that the wrapping number of this knot in the solid torus
$N(\gamma)$ is $\geq 2$ \cite{E2}). Note that if such a knot is not well wrapped then $S$ is in
fact meridionally compressible.

\begin{thm}\label{thm:type2} Let $\Gamma_2$ be a graph as above, and $k$ in $N(\Gamma_2)$ 
as above. 
Then $S=\partial N(\Gamma_2)$ is meridionally incompressible in $S^3-k$, 
and $k$ is a $(1,2)$--knot.
\end{thm}

\begin{proof} By construction $k$ is in a  $(1,2)$--position. The surface $S$ 
is meridionally incompressible in $N(\Gamma_2)-k$ by hypothesis. So it remains to prove 
that $S$ is incompressible in $S^3-\Int N(\Gamma_2)$.

Let $A=T_e -\Int N(\Gamma_2)$. This is a once-punctured annulus if $\gamma$ is nontrivial
in $T_e$, and it is a once-punctured disk plus a once-punctured torus if $\gamma$ is trivial in $T_e$.
Suppose $D$ is a compression disk for $S$. Look at the intersections between $D$ and $A$.
Let $\beta$ be a simple closed curve of intersection which is innermost in $D$. 
If $\beta$ is nontrivial in $T_e$, then it is a meridian of the solid torus 
$(T\times [e,1])\cup R_1$, say, which  implies that the arc $\alpha_2$ can be isotoped 
into $T_e$ disjoint from $\gamma$.
So suppose that $\beta$ is trivial in $T_e$ (but perhaps nontrivial in $A$).

If $\gamma$ is
nontrivial in $T_e$, then the curve $\beta$ is trivial in $A$, and it is easily removed.
If $\gamma$ is trivial in
$T_e$, then $\beta$ will be trivial in $A$, except if it is a curve concentric with $\gamma$, 
not contained in the disk bounded by $\gamma$. In this case, the disk bounded by $\beta$ in $D$ 
union the disk bounded by $\beta$ in $T_e$ bounds a $3$--ball which contains the upper or the lower 
arc of $\alpha$, and then as such arc has no local knots (for it has just one maximum or minimum), 
it can be pushed into $T_e$, which
contradicts our hypothesis. Suppose then that all simple closed curves of intersection 
between $A$ and $D$ have been removed.

Suppose now that $\beta$ is an arc of intersection between $D$ and $A$ which is outermost in $D$.
Then $\beta$ cuts off a disk $D'\subset D$, with $\partial D'=\beta\cup \delta$, 
where $\delta \subset S$. 
If $\beta$ is trivial in $A$, ie isotopic into a component of $\partial A$, 
then by cutting $D$ with the disk in $A$ determined by $\beta$ (or an innermost one), 
we get a compression disk with fewer intersections with $A$. Assume then that $\beta$ is 
nontrivial in $A$. Suppose first that the ends of $\beta$ lie on
$N(\gamma)$; in this case we can assume that $\delta$ is disjoint from $\partial N(\alpha)$
(for otherwise the interior of $\delta$ would intersect $A$). If $\gamma$ is
trivial in $T_e$, then $D'$ is a meridian disk of $(T\times [e,1])\cup R_1$, say, which implies
that the arc $\alpha_2$ can be isotoped into $T_e$. If $\gamma$ is nontrivial in $T_e$, then 
either $D'$ is a meridian disk of $(T\times [e,1])\cup R_1$, say, and $\gamma$ is a curve intersecting a
meridian of $(T\times [e,1])\cup R_1$ in one point, or $\partial D'$ determines a disk $E\subset T_e$,
which contains the point $T_e\cap \alpha$.
Again, in both cases this implies that the arc $\alpha_2$ can be isotoped into $T_e$.
Suppose now that both ends of $\beta$ lie on $N(\alpha)$. If the arc $\delta$ is 
isotopic into $\partial A$, then by isotoping $D$ we would get a compression
disk whose intersection with $T_e$ contains a simple closed curve. Otherwise, 
$\delta=\delta_1 \cup \delta_2 \cup \delta_3$, where $\delta_1$ and $\delta_3$ are arcs lying on 
$\partial N(\alpha)$, and $\delta_2$ is an arc lying on $\partial N(\gamma)$. Note that $\delta_2$
goes around $N(\gamma)$ just once. Let $E$ be a disk contained in $\partial N(\alpha)$,
such that $\partial E= \delta_1 \cup \delta_4 \cup \delta_3 \cup \delta_5$, where 
$\delta_4 \subset \partial N(\alpha)\cap \partial N(\gamma)$, and $\delta_5 \subset \partial A$.
Then $D'\cup E$ is an annulus in $(T\times [e,1])\cup R_1$, a boundary component of it
is a curve parallel to $\gamma$, the other component lies on $A$, and the arc $\alpha_2$
is an spanning arc of $D'\cup E$. This again shows that $\alpha_2$ can be isotoped into $T_e$.
Finally, if one endpoint of $\beta$ lies in $N(\gamma)$ and the other in $N(\alpha)$, then
$\alpha_2$ can be isotoped to lie on $T_e$. 
\end{proof}

Note that if in the $1$--complex $\Gamma_2$ one of the arcs $\alpha_1$ or $\alpha_2$ can  be isotoped
into an arc on the level torus $T_e$ with interior disjoint from $\gamma$, then the surface $S$
will be either compressible, or it can be isotoped to a surface of type 1, so that the knot $k$ remains
in a $(1,2)$--position.

\subsection{Surfaces of type 3}\label{sec:type3}

Let $e$ be a point on $I$, so that
$0 < e < 1$. Consider the torus $T_e = T\times \{e\}$.
Let $\gamma_1$ be a simple closed curve embedded in the
level torus $T_e$ and let $\gamma_2$ be an essential simple closed curve
embedded in $T_0$ which goes around $R_0$ at least once longitudinally. 
Let $\alpha$ be an arc contained in
$T\times I$, with endpoints in $\gamma_1$ and $\gamma_2$, so that it has just a local maximum at $T_1$.
Suppose that $\gamma_1$ is essential in $T_e$, or that it is 
inessential but bounds a disk in $T_e$ which intersects $\alpha$ in one point.
Note that the interior of $\alpha$ intersects $T_e$ in one point, so $\alpha$ 
is divided into a lower and an upper arc, say $\alpha_1$ and $\alpha_2$. 
Suppose that the arc $\alpha_2$ cannot be isotoped (in $(T\times [e,1])\cup R_1$),
keeping its endpoints fixed, into an arc on $T_e$ with interior disjoint from $\gamma_1$.

\begin{figure}[ht!]
\labellist
\small\hair 2pt
\pinlabel $T_1$ at 8 121
\pinlabel $T_0$ at 69 69
\pinlabel $\gamma_2$ at 59 51
\pinlabel $\gamma_1$ at 27 69
\pinlabel $\alpha$ at 19 102
\pinlabel $\gamma_1$ at 232 93
\pinlabel $\gamma_2$ at 169 72
\pinlabel $\alpha$ at 314 107
\endlabellist
\centering
\includegraphics{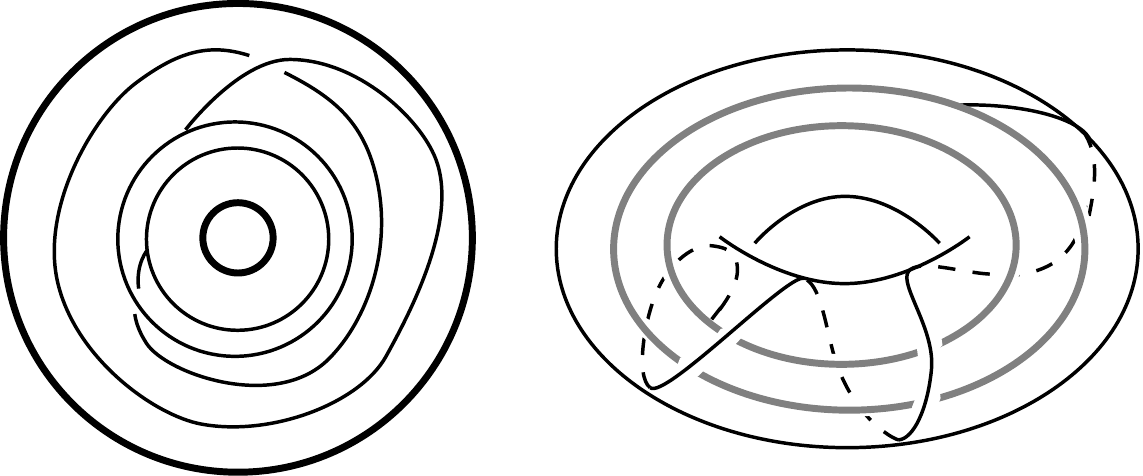}
\caption{}\label{fig:type3-4}
\end{figure}

Let $\Gamma_3$ be the $1$--complex consisting of
the union of  the curves $\gamma_1$, $\gamma_2$ and the arc $\alpha$. So
$\Gamma_3$ is a trivalent graph embedded in $S^3$. 
Let $N(\Gamma_3)$ be a regular neighborhood of
$\Gamma_3$. This is a genus $2$ handlebody. We can assume 
that $N(\Gamma_3)$
is the union of the solid tori $N(\gamma_1)$ and $N(\gamma_2)$, joined
by the $1$--handle $N(\alpha)$, so that these intersect in two disks, $E_1$ 
and $E_2$, where say $E_1$ is at level $T\times\{e+\epsilon\}$ and $E_2$ at level 
$T\times\{\epsilon\}$, for some small $\epsilon > 0$. 
Let $S=\partial N(\Gamma_3)$, we say that $S$ is a surface of type 3.
A surface of type 3 is shown schematically in \fullref{fig:type3-4}, and an explicit 
example is also given.

Let $k$ be a knot embedded in $N(\Gamma_3)$ so that $k$ intersects each of 
$E_1$ and $E_2$ in two points, $k\cap N(\alpha)$ consists of two arcs, each
with just a local maximum, and $k\cap N(\gamma_i)$, $i=1,2$, consists of one arc with
just a local minimum. Suppose also that
$S$ is meridionally incompressible in $N(\Gamma_3)-k$. To get that it suffices 
to ask that $k$ is well wrapped in $N(\Gamma_3)$, that is, consider
the arc of $k_i$ contained in $N(\gamma_i)$, join its endpoints lying in $E_i$ with an arc in $E_i$, 
push the resulting knot 
into $N(\gamma_i)$, and assume that the wrapping number of such a knot in $N(\gamma_i)$ is $\geq 2$.

\begin{thm}\label{thm:type3} Let $\Gamma_3$ be a graph as above, and $k$ in $N(\Gamma_3)$ 
as above. 
Then $S=\partial N(\Gamma_3)$ is meridionally incompressible in $S^3-k$, 
and $k$ is a $(1,2)$--knot.
\end{thm}

\begin{proof} By construction $k$ is in a  $(1,2)$--position. The surface $S$ 
is meridionally incompressible in $N(\Gamma_3)-k$ by hypothesis. 
So it remains to prove that $S$ is incompressible in $S^3-\Int N(\Gamma_3)$.

The proof is an innermost disk/outermost arc argument, looking at the intersections 
of a compression disk $D$ with the surface $T_e-\Int (N(\gamma_1)\cup N(\alpha))$.
\end{proof}

Note that if the curves $\gamma_1$ and $\gamma_2$ have the same slope, then $N(\Gamma_3)$ can be isotoped
so that both curves lie on the torus $T_0$. Note that if the arc 
$\alpha_2$ can  be isotoped
into an arc on the level torus $T_e$ with interior disjoint from $\gamma_1$, then the surface $S$
will be either compressible, or it can be isotoped to a surface of type 1, so that the knot $k$ remains
in a $(1,2)$--position.

\subsection{Surfaces of type 4}\label{sec:type4}

Let $e$, $\epsilon$ be points on $I$ so that
$0 < \epsilon < e < 1$. Consider the tori $T_e = T\times \{e\}$
and $F=T\times \{\epsilon\}$.
Let $\gamma_1$ be a simple closed curve embedded in the
level torus $T_e$ and let $\alpha$ be an arc contained in
$T\times I$, with endpoints in $\gamma_1$ and $F$, so that it has just a local maximum at $T_1$.
Suppose that $\gamma_1$ is essential in $T_e$, or that it is 
inessential but bounds a disk in $T_e$ which intersects $\alpha$ in one point.
Note that the interior of $\alpha$ intersects $T_e$ in one point, 
so $\alpha$ 
is divided into a lower and an upper arc, say $\alpha_1$ and $\alpha_2$. 
Suppose that the arc $\alpha_2$ cannot be isotoped (in $(T\times [e,1])\cup R_1$),
keeping its endpoints fixed, into an arc on $T_e$ with interior disjoint from $\gamma_1$.

The torus $F$ bounds a solid torus $F'=R_0\cup (T\times[0,\epsilon])$. Consider the union $H_4=N(\gamma_1)\cup
N(\alpha)\cup F'$. This is a genus 2 handlebody. This can be seen as the solid tori 
$N(\gamma_1)$ and $F'$
joined by the $1$--handle $N(\alpha)$, so that these 
intersect in two disks, 
$E_1$ and $E_2$, where say $E_1$ is at level $T\times\{e+\epsilon\}$,
and $E_2$ is at level $T\times\{\epsilon\}$. Let $S=\partial H_4$; we say that $S$ is a surface of type 4.
For an example of a surface of type 4, look at \fullref{fig:type3-4}, thinking of $\gamma_2$ as a fat solid
torus engulfing all of $R_0$.

Let $k$ be a knot embedded in $H_4$ so that $k$ intersects each of $E_1$ 
and $E_2$ in two points, $k\cap N(\alpha)$ consists of two arcs, each with
just a local maximum, $k\cap N(\gamma_1)$ consists of one arc, with just a
local minimum, and $k\cap F'$ consists also of one arc, with just a
local minimum. Suppose also that $S=\partial H_4$ is meridionally incompressible in $H_4$. 
To get that it suffices to ask that $k$  is well wrapped in $H_4$.

\begin{thm}\label{thm:type4}
Let $H_4$ be as above, and $k$ in $H_4$ as above. 
Then $S=\partial H_4$ is meridionally incompressible in $S^3-k$, and 
$k$ is a $(1,2)$--knot.
\end{thm}

\begin{proof} By construction $k$ is in a $(1,2)$--position, and by hypothesis $S$ 
is meridionally incompressible in $H_4-k$. So it remains to prove that $S$ is
incompressible in the complement $S^3-\Int H_4$. Such a proof is again an 
innermost disk/outermost arc argument, looking at the intersections 
between a compression disk $D$ and the surface $T_e-\Int N(\gamma_1\cup \alpha)$.
\end{proof}

Note that a surface of type 4 can be isotoped to look like a surface of type 3, 
where the curve $\gamma_2$ 
for the new surface of type 3 will be longitudinal. But if this isotopy is done then the  knot $k$ 
constructed for the surface of type 4 may not be
in a $(1,2)$--position. But a surface of type 3, where the curve $\gamma_2$ is longitudinal, will
be in fact a surface of type 4. Note also that if the arc $\alpha_2$ can  be isotoped
into an arc on the level torus $T_e$ with interior disjoint from $\gamma_1$, then the surface $S$
will be compressible.

\subsection{Surfaces of type 5}\label{sec:type5}

Consider a sphere $\Sigma$ that consists of two meridian disks in $R_1$, say $D_1$ and $D_2$, 
two vertical annuli $A_1$ and $A_2$ in 
$T\times [\epsilon,1]$, $0< \epsilon <1$, and an annulus $A_3$ in the
level torus $T\times \{\epsilon\}$. 
Let $B$ be a $3$--ball bounded by $\Sigma$ in $S^3$, say the one which does not contain the point at infinity. 

Assume first that the solid torus $R_0$ is not contained in $B$. 
Let $\gamma$ be a level simple closed curve lying in a level torus $T_e = T\times \{e\}$,
$0< \epsilon < e < 1$, and which lies inside the $3$--ball $B$. Let $\alpha_1$ be a vertical arc in $B$ with
one endpoint in $A_1$ at a level above $T_e$, and the other endpoint in $\gamma$, and let $\alpha_2$ be a
vertical arc in $B$ with one endpoint in $A_2$ and the other in $A_2$ or in $A_1$. If $\gamma$ is a trivial
curve in the level torus $T_e$, assume that $\alpha_2$ intersects in one point the level disk $E$ bounded by
$\gamma$. Assume that there is no disk $D$ in $B$ with $\partial D=\alpha_2\cup \delta$, $\delta\subset \Sigma$,
and $D\cap (\alpha_1\cup\gamma)=\emptyset$, that is, the arc $\alpha_2$ is not isotopic to an arc in $\Sigma$. 
Suppose also that the arc $\alpha_1$ cannot be isotoped, keeping its endpoints in $\Sigma$ and $\gamma$, so that
$\alpha_1$ lies in the level torus $T_e$, and the arc $\alpha_2$ remains being a vertical arc. It is not 
difficult to construct examples satisfying these conditions. 

\begin{figure}[ht!]
\labellist
\small\hair 2pt
\pinlabel $D_1$ at 110 141
\pinlabel $D_2$ at 111 9
\pinlabel $\gamma$ at 26 69
\pinlabel $\alpha_1$ at 76 108
\pinlabel $\alpha_2$ at 60 40
\pinlabel $A_1$ at 173 110
\pinlabel $A_2$ at 173 39 
\pinlabel $A_3$ at 6 65
\endlabellist
\centering
\includegraphics{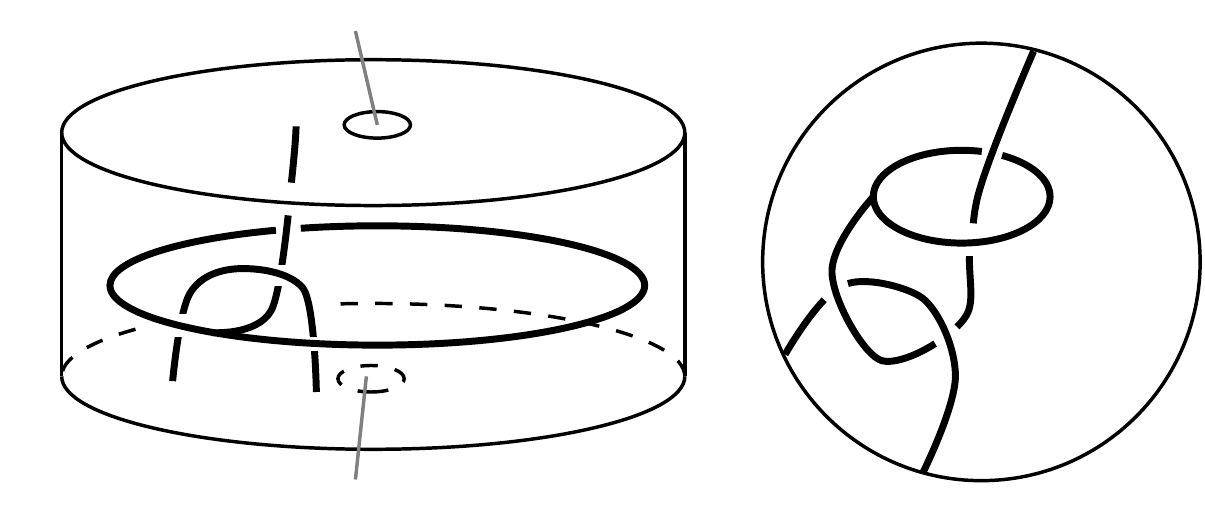}
\caption{}\label{fig:type5}
\end{figure}

Let $B_1$ be the complementary ball in $S^3$ bounded by $\Sigma$. 
Let $H_5=B_1\cup N(\alpha_1)\cup N(\alpha_2)\cup N(\gamma)$.
This is a genus 2 handlebody. First, $B_1\cup N(\alpha_2)$ is a solid torus, formed by the $3$--ball $B_1$ and
the $1$--handle $N(\alpha_2)$, where $B_1\cap N(\alpha_2)$ consists of two vertical disks, say $E_1$ and
$E_2$, where $E_1$ is at an upper level. So $H_5$ can be seen as the solid tori $B_1\cup N(\alpha_2)$ and
$N(\gamma)$ joined by the $1$--handle $N(\alpha_1)$, where $B_1\cap N(\alpha_1)$ consists of a vertical disk
$E_3$, and $N(\alpha_1)\cap N(\gamma)$ is a level disk $E_4$, lying in a level $T\times \{e+\delta\}$,
for some small $\delta > 0$. Let $S=\partial H_5$. We say that $S$ is a surface of type 5. Look at 
\fullref{fig:type5}, left, for an example of a surface of type 5.

Let $\Sigma$ and $B$ be as above, but suppose now that the solid torus $R_0$ is contained in the $3$--ball $B$.
Let $F=T\times \{\epsilon_1\}$, $0<\epsilon_1 < \epsilon$, and let $F'=R_0\cup(T\times [0,\epsilon_1])$.
So the solid torus $F'$ is contained in $B$. Let $\alpha_1$ be a vertical arc in $B$ with
one endpoint in $A_1$ and the other endpoint in $F$, and let $\alpha_2$ be a
vertical arc in $B$ with one endpoint in $A_2$ and the other in $A_2$ or in $A_1$.
As before, assume that the arc $\alpha_2$ is not isotopic to an arc in $\Sigma$.
Again, let $B_1$ be the complementary ball in $S^3$ bounded by $\Sigma$, and let $H_5$ be the genus 2 handlebody 
$H_5=B_1\cup N(\alpha_1)\cup N(\alpha_2)\cup F'$. Define disks $E_1$, $E_2$, $E_3$ and $E_4$ as above.
Let $S=\partial H_5$, we also say that $S$ is a surface of type 5.

Let $k$ be a knot in $H_5$, intersecting in two points each of the disks $E_i$, $i=1,2,3,4$, and so that 
$k\cap N(\alpha_2)$ consists of two vertical arcs, $k\cap N(\alpha_1)$ consists of two vertical arcs,
$k\cap N(\gamma)$ (or $k \cap F'$) consists of one arc with a single local minimum and which is well wrapped in
$N(\gamma)$ ($F'$), and $k\cap B_1$ consists of three arcs, two of them with a single local
maximum and with endpoints in $E_1\cup E_3$, the other with a single local minimum
and with endpoints in $E_2$. Suppose also that none of these arcs is isotopic in $B_1$, keeping its endpoints 
and the other arcs fixed, to an arc lying on some $E_i$, $i=1,2,3$. 

Note that the two constructions of surfaces of type 5 produce surfaces which are isotopic in $S^3$, but if
such an isotopy is performed transforming one surface then the corresponding knot $k$ may not longer be in a
$(1,2)$--position.

\begin{thm}\label{thm:type5}
Let $S$ and $k$ be as above. $S$ is meridionally incompressible in 
$S^3-k$, and $k$ is a $(1,2)$--knot. 
\end{thm}

\begin{proof} 
By construction $k$ is in a
$(1,2)$--position. We have to show that
$S$ is incompressible in 
$B- \Int H_5$ and that $S$ is meridionally incompressible in $H_5-k$. It is not difficult to prove that $S$ 
is meridionally incompressible in $H_5-k$. To do that  
consider the disks $E_1$, $E_2$, $E_3$ and $E_4$; these are disks which intersects $k$ in two points. 
Look at the intersections between a compression disk $D$ and $E_1\cup E_2\cup E_3\cup E_4$. Using the 
hypothesis on $k$, we conclude that $D$ and $E_1\cup E_2\cup E_3\cup E_4$ can be made disjoint, which then 
implies that $D$ cannot exist. 

Suppose that the solid torus $R_0$ is not contained in the $3$--ball $B$, the proof for the remaining case is
similar. Note that if the lower endpoint of the arc $\alpha_2$ lies in the annulus
$A_1$, and it is at a level below
$T_e$, then it can be isotoped, going through the annulus $A_3$, so that
both of its endpoints lie in $A_2$. This isotopy can be performed, moving $H_5$ and $k$, but so 
that $k$ remains in a $(1,2)$--position. If $\gamma$ is a nontrivial curve in the level torus $T_e$, and both
endpoints of $\alpha_2$ lie in $A_2$, then $\gamma$ can be isotoped so that it lies at a level below the lower
endpoint of $\alpha_2$. Assume these isotopies have been performed, if possible.

Let $E=T_e \cap (B-\Int H_5)$. If $\gamma$ is a nontrivial curve in the level torus $T_e$, 
then $E$ consist of two annuli. And if $\gamma$ is a trivial curve in the level torus then $E$ consists of a
punctured annulus and a punctured disk. The proof is now an innermost disk/outermost arc argument, looking at the
intersections between $E$ and a compression disk.
\end{proof}

Note that if the arcs $\alpha_1$ and $\alpha_2$ do not satisfy the required hypothesis, then the surface $S$
will be compressible.

\subsection{Surfaces of type 6}\label{sec:type6}

Let $\gamma$ be a knot in $T\times I$ in a $(1,1)$--position, so that it has
a local maximum at a level just below $T_1$, and a local minimum at a level just above $T_0$.
Assume that $\gamma$ is not isotopic in $T\times I$ to a meridian or a longitude of
a level torus. 
Let $N(\gamma)$ be a neighborhood of $\gamma$. Let $\alpha_1$ be a trivial curve, 
in a level torus $T\times \{ e\}$, $0< e < 1$, which bounds a
level disk $E$ such that $E \subset \Int N(\gamma)$. Let $\alpha_2$ be an arc contained in
$N(\gamma)$ with an endpoint in $\partial N(\gamma)$, lying at a level $f$,
with $0 < f < e$, and the other point in $\alpha_1$. Suppose also that $\alpha_2$ has
a single local maximum in $T\times I$, that 
$\alpha_2$ intersects in one point the disk $E$ bounded by $\alpha_1$, and that
$\alpha_1\cup \alpha_2$ intersects each meridian of $N(\gamma)$, that is, we have 
something like in \fullref{fig:type6}. Assume that $N(\alpha_1)\cap N(\alpha_2)$ consists of a level disk
$E_1$ lying in a level torus $T\times \{e+\epsilon\}$, for some $\epsilon > 0$, and that
$N(\alpha_2)\cap \partial N(\gamma)$ is a disk $E_2$ lying in the level torus $T\times \{f\}$.
Let $M_6=N(\gamma)-\Int N(\alpha_1)\cup N(\alpha_2)$. Note that $S=\partial M_6$ is a 
genus 2 surface. We say that $S$ is a surface of type 6.

Let $k$ be a knot contained in $S^3-M_6$, which intersects each of $E_1$, $E_2$ in two points, and 
so that $k\cap N(\alpha_2)$ consists of two arcs each having a 
single local maximum, $k\cap N(\alpha_1)$ consists of one arc, well wrapped in 
$N(\alpha_1)$, and having just one local minimum, and $k\cap (S^3-\Int N(\gamma))$ is
one arc, with a single local minimum and which goes 
around a longitude of $R_0$ at
least once. If $\gamma$ is a trivial knot assume further that $k$ is well wrapped in
$S^ 3-\Int N(\gamma)$. See \fullref{fig:type6}.

\begin{thm}\label{thm:type6}
Let $S$ and $k$ be as above. $S$ is meridionally incompressible in $S^3-k$, 
and $k$ is a $(1,2)$--knot. Furthermore, if $\gamma$ is a nontrivial knot, then $S$ does not 
bound a handlebody in $S^3$.
\end{thm}

\begin{proof} We have to show that $S$ is incompressible in $M_6$ and that $S$ is meridionally 
incompressible in $(S^3-M_6)-k$. It is not difficult to prove that $S$ is meridionally 
incompressible in $(S^3-M_6)-k$. To do that look at the intersections between a compression disk $D$ 
and $E_1\cup E_2$. 
Using the facts that $k$ is well wrapped in $N(\alpha_1\cup \alpha_2)$ and that $k$ goes at least once
longitudinally around $R_0$, or that $k$ is well wrapped in $S^ 3-\Int N(\gamma)$, we conclude that 
$D$ and $E_1\cup E_2$ can be made disjoint, which then implies
that $D$ cannot exist. 

\begin{figure}[ht!]
\labellist
\small\hair 2pt
\pinlabel $N(\gamma)$ at 83 157
\pinlabel $\alpha_1$ at 43 52
\pinlabel $\alpha_2$ at 31 134
\pinlabel $k$ at 325 22
\endlabellist
\centering
\includegraphics[scale=.9]{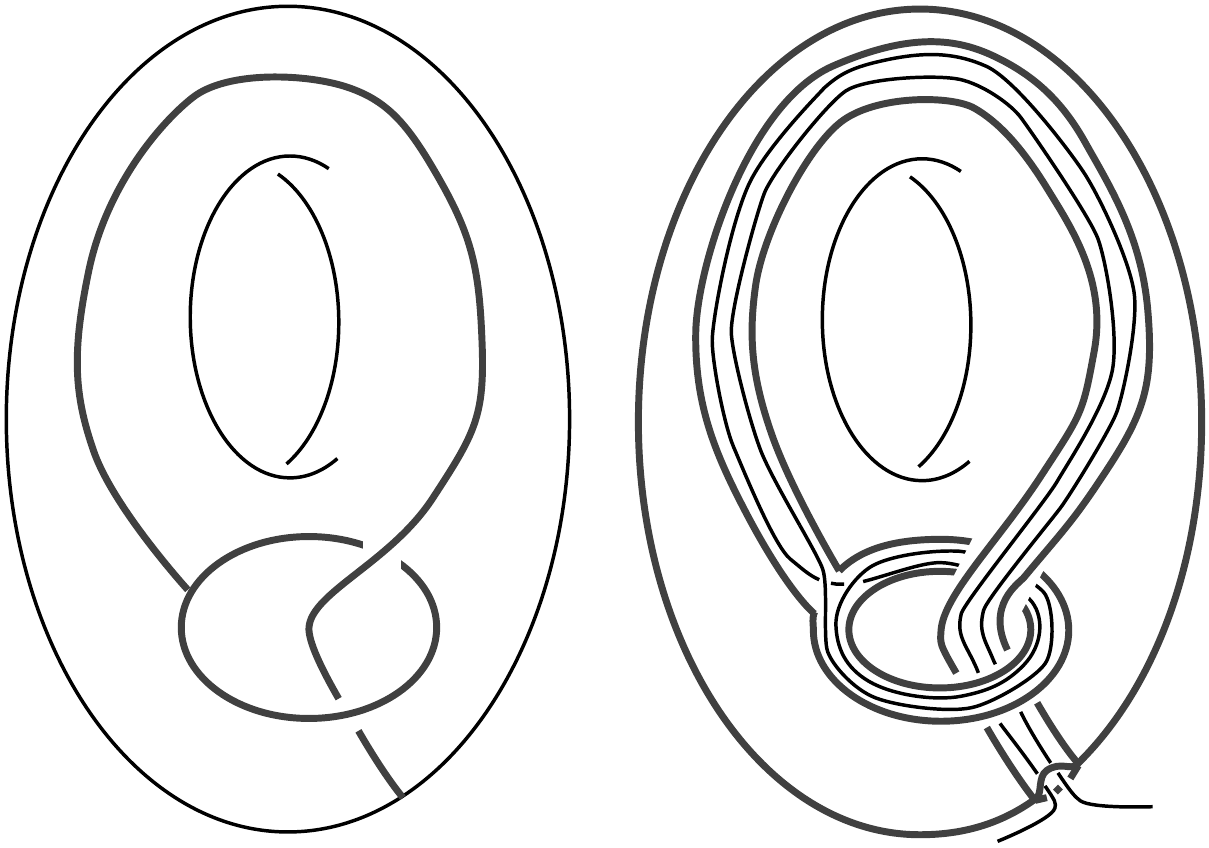}
\caption{}\label{fig:type6}
\end{figure}

It remains to prove that $S$ is incompressible in $M_6$. Note that there are two nonseparating 
annuli properly embedded in $M_6$, say $A_1$ and $A_2$, so that the boundary of $A_1$ consists of a
longitude of $\alpha_1$ and the boundary of a cocore of the $1$--handle 
$N(\alpha_2)$, and the boundary of $A_2$ consists of a meridian of
$\gamma$ and the boundary of a cocore of the $1$--handle 
$N(\alpha_2)$. Note that anyone of the boundary components of $A_1$ and $A_2$ is a nontrivial curve in $M_6$, 
because $\alpha_1\cup \alpha_2$ is not contained in a $3$--ball inside
$N(\gamma)$. In particular this shows that the annuli $A_1$ and $A_2$
are incompressible. Note also that
$A_1\cup A_2$ does not separate $M_6$. Suppose that $D$ is a compression disk for
$S$ and look at the intersections between $D$ and $A_1\cup A_2$. Simple closed curves of 
intersection are easily
removed, for these have to be trivial in the annuli. So the intersection consists only of arcs.
Let $\beta$ be an outermost arc of intersection in $D$, so it cuts off a disk $D'\subset D$, with
$\partial D'=\beta\cup \delta$. If the endpoints of $\beta$ lie on different boundary components of 
$A_1$ or $A_2$,
then by inspection we see that there cannot be an arc $\delta$ in $S$, with interior disjoint from 
the annuli, joining these two points. So if such $\beta$ exists, it must have endpoints in the same 
boundary component of one of the annuli, so it bounds a disk $D''$ in the annulus, and by cutting $D$ with
$D''$  (or with one outermost disk contained in $D''$), we get a disk with fewer intersections with the
annuli. Note also that $D$ cannot be disjoint from the annuli, for otherwise $\alpha_1\cup \alpha_2$
will be contained in a $3$--ball inside $N(\gamma)$. 

Finally note that if $\gamma$ is a nontrivial knot, then $S$ does not bound a handlebody in $S^3$, 
for in one it side bounds the disk sum of $S^3-\Int N(\gamma)$ with $N(\alpha_1) \cup N(\alpha_2)$, and the other
side it bounds the manifold $M_6$, which is not a handlebody for it has incompressible boundary.
\end{proof}

Note that a knot $k$ whose complement has a surface of type 6 will in general also have a surface of type 4. To see that, consider the union of the curve $\alpha_1$, the arc $\alpha_2$ and $R_0$ 
(where the arc $\alpha_2$ has been prolonged to touch $R_0$). So, let $H_4=N(\alpha_1\cup \alpha_2\cup
R_0)$. Note that $S'=\partial H_4$ is a surface of type 4, and if $H_4$ is thin then
$S\cap S' =\emptyset$. The surface $S'$ will be in fact meridionally incompressible in $S^3-k$, 
except if $k$ goes around $R_0$ exactly once longitudinally. 

\subsection{Surfaces of type 7}\label{sec:type7}

\begin{figure}[ht!]
\labellist
\small\hair 2pt
\pinlabel $N_1$ at 71 62
\pinlabel $N_2$ at 71 17
\pinlabel $\alpha_1$ at 146 78
\pinlabel $\alpha_2$ at 229 97
\pinlabel $A_2$ at 24 20
\pinlabel $A_2$ at 300 20
\endlabellist
\centering
\includegraphics{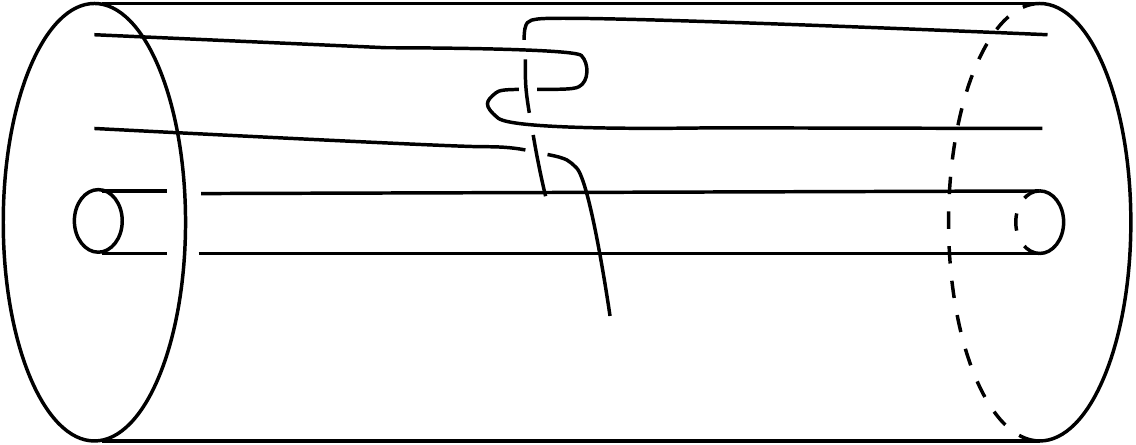}
\caption{}\label{fig:type7}
\end{figure}

Let $\gamma$ be a simple closed curve in $T\times\{1/2\}$ of slope 
$p/q$, $\vert p\vert \geq 1$ 
(in the usual coordinates for the solid torus $R_0\cup (T\times\{1/2\})$. Let $N_1=N_1(\gamma)$ and 
$N_2=N_2(\gamma)$ be two regular neighborhoods of $\gamma$, with $N_1 \subset N_2$, and 
$N_1 \subset T\times [1/4,3/4]$, $N_2 \subset T\times [1/8,15/16]$. Let $\alpha$ be an arc contained
in $N_2 -\Int N_1$, connecting $\partial N_1$ with $\partial N_2$, and so that $\alpha$ has just
one local maximum. Suppose also that $\alpha$ cannot be isotoped, keeping its endpoints in
$\partial N_1$ and $\partial N_2$, to an arc lying on a level torus $T\times \{y\}$.
The arc $\alpha$ can be isotoped so that it looks like the union of the arcs $\alpha_1\cup \alpha_2$, 
where $\alpha_1$ is an arc of the form $\{x\}\times [1/2,7/8]$, going from $\partial N_1$ to the local
maximum, and $\alpha_2$ is a descending arc, which wraps around $N_2$ and around $\alpha_1$, until it
finishes at a point in $\partial N_2$. See \fullref{fig:type7} for an example. By isotoping and sliding the
arc $\alpha$, and maintaining it with a single maximum, we can assume that its endpoints lie in the same 
level, say at level $1/2$. We can connect its endpoints with an arc $\beta$ contained in an annulus $A$ 
(which is one of the components of  
$(N_2 -\Int N_1)\cap (T\times \{1/2\})$), and in fact there are many of such arcs. Assuming that the arc
$\alpha$ cannot be isotoped to be level, is equivalent to assuming that the knot 
$\ell = \alpha\cup \beta$ is never trivial in $N_2 -\Int N_1$, for any of the choices of $\beta$. We can think 
of the knot $\ell$ as lying in $A\times I$, so that $\ell$ has just a maximum in
$A\times I$. By embedding $A\times I$ as a standard solid torus in $S^3$, with $\partial A\times \{0\}$ 
being a preferred longitude of such a solid torus, we get a $2$--bridge link, formed by $\ell$ and a meridian
of the solid torus 
$A\times I$. The assumption on the arcs is equivalent to asking that the corresponding $2$--bridge 
link is never a split link.

We can assume that a neighborhood of $\alpha$, $N(\alpha)$, is contained in $N_2-\Int N_1$, and that
$N_1\cap N(\alpha)$ consists of a disk $E_1$ and that $N_2\cap N(\alpha)$ is a disk $E_2$.
Let $M_7=N_2-\Int (N_1\cup N(\alpha))$, and let $S=\partial M_7$. We say that $S$ is
a surface of type 7. Let $k$ be a knot contained in $S^3 -M_7$, which intersects each of $E_1$, $E_2$ in
two points, and so that $k\cap N(\alpha)$ consists of two
arcs, each with a single local maximum, and each with an endpoint in $E_1$ and the other in $E_2$. $k\cap N_1$
is an arc with endpoints in $E_1$, with a single local minimum, and which is well wrapped in $N_1$, and 
$k \cap(S^3-N_2)$ is an arc with endpoints in $E_2$, with a single local minimum, and which goes around 
$R_0$ at least once longitudinally. If $\gamma$ is a trivial knot assume further that $k$ is well wrapped
in $S^3-\Int N_2$.

\begin{thm}\label{thm:type7}
Let $S$ and $k$ be as above. $S$ is meridionally incompressible in 
$S^3-k$, 
and $k$ is a $(1,2)$--knot. Furthermore, if $\gamma$ is a nontrivial knot, then $S$ does not bound a 
handlebody in $S^3$.
\end{thm}

\begin{proof} 
It is not difficult to prove that $S$ is meridionally 
incompressible in $(S^3-M_7)-k$. Look at the intersections between a compression disk $D$ and the 
disks $E_1$ and $E_2$. 
Using the hypothesis on $k$, we conclude that $D$ and $E_1$, $E_2$ can be made disjoint, 
which then implies that $D$ cannot exist.

It remains to prove that $S$ is incompressible in $M_7$. $M_7\cap (T\times \{ 1/2 \})$ consists of two
annuli, and $\alpha$ must be disjoint from one of these annuli, for otherwise it will contain more than
one local maxima. So let $A_1$ be one of such annuli, and suppose that 
$A_1$ is disjoint from $\alpha$ (the other annulus was denoted 
before by $A$).
$A_1$ is a nonseparating annulus in $M_7$. Let $A_2$ be an annulus in $N_2-\Int N_1$ consisting of a
meridian disk of $N_2$ minus a meridian disk of $N_1$.
$A_1$ and $A_2$ intersect in a single arc which is essential in both annuli. The arc $\alpha$ must 
intersect $A_2$, for otherwise it will be contained in a $3$--ball and it would be isotopic to a level arc. 
So suppose that 
$\alpha$ and $A_2$ intersect transversely and that this intersection is minimal. $A_2\cap M_7$
is then a punctured annulus, which we call also $A_2$.

Suppose that $S$ is compressible in $M_7$, and let
$D$ be a compression disk. Look at the intersections between
$D$ and $A_1\cup A_2$. Simple closed curves of intersection between $D$ and $A_1$ are easily removed, 
for no such curve can be essential in $A_1$. For the same reason, if there is a simple closed curve of
intersection between $D$ and $A_2$, it must bound a disk in $A_2$; if the disk intersects $\alpha$, 
then an isotopy reduces the number of  points of intersection 
between $A_2$ and $\alpha$, otherwise such an isotopy reduces the number of
intersection curves between $D$ and $A_2$. So the
intersection consists only of arcs. By isotoping
$D$, we can assume that it is disjoint from the arc of intersection between $A_1$ and $A_2$.
Let $\beta$ be an outermost arc of intersection between $D$ and $A_1\cup A_2$, and suppose that $\beta$ 
lies on $A_1$. If $\beta$ is inessential in $A_1$ then it is easily removed. 
So suppose $\beta$ has endpoints on different components of $A_1$. The arc
$\beta$ cuts off a disk $D'\subset D$, with $\partial D=\beta\cup \delta$. But the arc $\delta$ must 
pass through $N(\alpha)$, for otherwise cannot connect points on different components of $\partial A_1$.
Then $\delta$ must intersect $A_2$, which contradicts the fact that $\beta$ is outermost. So any 
outermost arc of intersection $\beta$ must lie on
$A_2$. In any of the possible cases, the disk $D'$ determined by $\beta$ can be used to isotope $\alpha$,
reducing the number of points of intersection of $\alpha$ with $A_2$. So the disk $D$ must be disjoint from 
$A_1$ and $A_2$. Then it is not difficult to see that this is not possible.

Finally, if $\gamma$ is a nontrivial knot, then $S$ does not bound a handlebody in $S^3$, for in one side 
it bounds the disk sum of $S^3-\Int N_2$ with $N(\alpha) \cup N_1$, and in the other side it bounds $M_7$, 
which is not a handlebody for it has incompressible boundary.
\end{proof}

Note that if the curve $\gamma$ is a 
curve of slope $1/q$ on $T\times \{1/2\}$, then $S$ will be
isotopic to a surface of type 4.
Note that a knot $k$ whose complement has a surface of type 7 will in general also have a surface of type 4.
To see that, consider the union of the solid torus $N_1$, the arc $\alpha$ and $R_0$ (where the arc
$\alpha$ has been prolonged to touch $R_0$). So, let $H_4=N_1 \cup N(\alpha)\cup R_0$.
Note that $S'=\partial H_4$ is a surface of type 4, and that if $H_4$ is thin enough then
$S\cap S' =\emptyset$. The surface $S'$ will be in fact meridionally incompressible in $S^3-k$, 
except if $k$ goes around $R_0$ exactly once longitudinally.

\subsection{Surfaces of type 8}\label{sec:type8}

Let $R$ be a torus in $S^3$ constructed as follows. Let $A_1$, $A_2$, $\dots$, $A_n$ be $n$ annuli properly 
embedded in $R_1$, all with slope $p/q$,  $\vert p\vert \geq 2$, $\vert q\vert \geq 2$ (in the usual 
coordinates of the solid torus $R_0\cup (T\times [0,1])$). Suppose the annuli are nested, and say, 
$A_1$ is the innermost one. Let $\smash{B_1^1}$, $\smash{B_1^2}$, $\dots$, $\smash{B_n^1}$, $\smash{B_n^2}$ be $2n$ vertical annuli 
in $T\times I$, so that $\partial A_i \subset \partial (\smash{B_i^1\cup B_i^2})$, $1\leq i\leq n$. 
Let $C_1$, $\dots$, $C_n$ be $n$ annuli properly embedded in $R_0$, which are nested, whose boundaries 
coincide with the boundaries of the $\smash{B_i^j}$'s, so that $C_1$ is the innermost annulus and 
$\partial C_1\subset \partial (\smash{B_1^2\cup B_2^2})$. Note that the union of the $A_i$'s, the $\smash{B_i^j}$'s 
and the $C_i$'s is an embedded torus, denoted by $R$. 

In the special case when $n=1$, assume that the annulus $C_1$ is chosen so that the torus $R$
does not bound a solid torus contained in $T\times I$. In the special case when $n=2$, it is enough to assume
that $\vert p\vert \geq 2$, $\vert q\vert \geq 1$. In this case we can take the annuli $C_1$ and $C_2$ to be
nested or non-nested.

Note that $R$ is a standard torus in $S^3$, except in the case when $n=2$, $\vert q\vert\geq 2$, and the annuli
$C_1$, $C_2$ are non-nested. In that case $R$ is isotopic to the boundary of a regular
neighborhood of the $(p,q)$--torus knot.

Let $\alpha$ be an arc in $T\times I$, with one endpoint in $\smash{B_1^1}$, the other in $B_1^2$, with 
interior disjoint from 
the $\smash{B_i^j}$'s, and so that it has a single local maximum in $T\times I$.
Let $E$ be the annulus in $T_0$, with $\partial E\subset \partial (\smash{B_1^1\cup B_1^2})$, and whose 
interior is disjoint 
from the $\smash{B_i^j}$'s. So $A_1\cup \smash{B_1^1\cup B_1^2}\cup E$ bounds a solid torus $P$ which contains the arc 
$\alpha$.  By sliding the arc $\alpha$, we can assume that its endpoints lie in the same level, say at level
$e$, $0 < e < 1$.  We assume that the arc $\alpha$ cannot be isotoped, keeping its endpoints fixed, to a level
arc lying in $T\times \{e\}$. 

\begin{figure}[ht!]
\labellist
\small\hair 2pt
\pinlabel $T_1$ at 5 47
\pinlabel $T_0$ at 82 29
\pinlabel $T_1$ at 194 47
\pinlabel $T_0$ at 271 29
\pinlabel $\alpha$ at 71 111
\pinlabel $\alpha$ at 250 115
\pinlabel $\gamma$ at 267 99
\endlabellist
\centering
\includegraphics{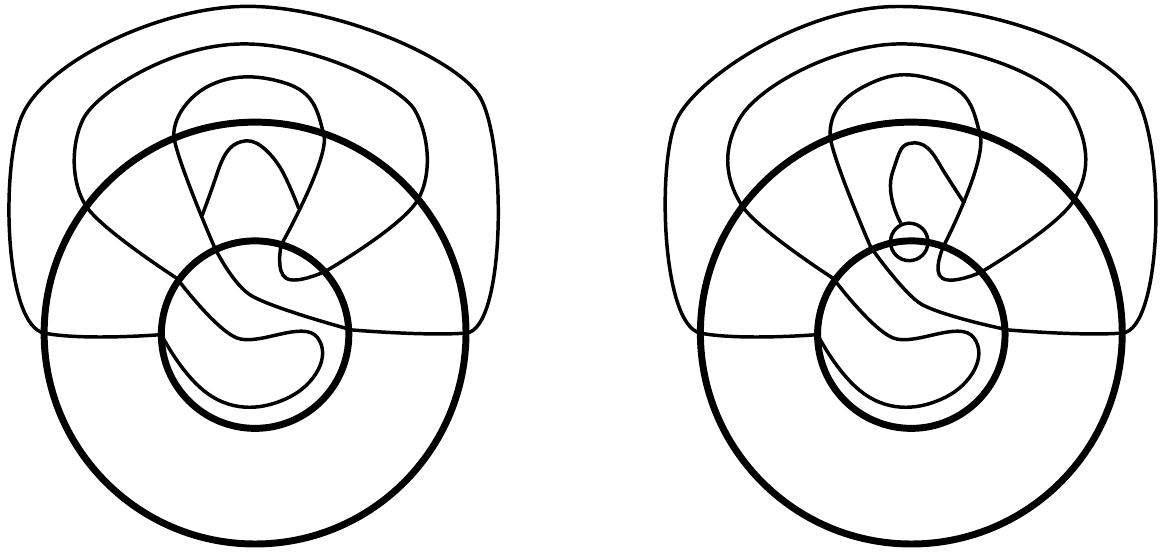}
\caption{}\label{fig:type8-9}
\end{figure}

Let $R'$ be the solid torus bounded by $R$ which does not contain the arc $\alpha$.
Let $H_8=R'\cup N(\alpha)$. This is a genus 2 handlebody; it can be seen as the solid torus $R'$ 
union the $1$--handle $N(\alpha)$, where $R'\cap N(\alpha)=(B_1^1\cup B_1^2)\cap N(\alpha)$ consists of two
disks, say $E_1$ and $E_2$. Let $S=\partial H_8$. We say that $S$ is a surface of type 8. See 
\mbox{\fullref{fig:type8-9}}, left.

Let $k$ be a knot in $H_8$, intersecting in two points each of $E_1$, $E_2$, and so that $k\cap N(\alpha)$ 
consists of two arcs, each with a single local maximum, and $k\cap R'$ consists of two arcs, each with a 
single local minimum,  and each going at least once longitudinally  around $R'$, ie none of the arcs can be
isotoped into an arc lying in $E_1$ or $E_2$. 

\begin{thm}\label{thm:type8}
Let $S$ and $k$ be as above. $S$ is meridionally incompressible 
in $S^3-k$, 
and $k$ is a $(1,2)$--knot. 
\end{thm}

\begin{proof} By construction $k$ is in a $(1,2)$--position. We have to show that $S$ is incompressible in 
$S^3- \Int H_8$ and that $S$ is meridionally incompressible in $H_8-k$. It is not difficult to prove that $S$ is
meridionally incompressible in $H_8-k$. To do that consider the disks $E_1$, $E_2$, and look at the 
intersections between a compression disk $D$ and $E_1\cup E_2$. Using the hypothesis on $k$, we conclude that $D$
and $E_1\cup E_2$ can be made disjoint, which then implies that $D$ cannot exist. 

Suppose now that $S$ is compressible in $S^3- \Int H_8$, and let $D$ be a compression disk. Look at the 
intersections  between $D$ and the annulus $E$ defined above. Simple closed curves of intersection are easily
removed. Let $\beta$ be an outermost arc of intersection in $D$, which cuts a disk $D'\subset D$, and say 
$\partial D'= \beta \cup \delta$. If $\delta$ is trivial in $E$, we cut $D$ with the disk in $E$ determined by 
$\delta$ (or with an innermost one), getting a compression disk with fewer intersections with $E$.
If $\delta$ is essential in $E$, there are two cases. If $D'$ is contained in the solid torus $P$, 
this will imply that the arc $\alpha$ can be isotoped to be level. If $D'$ is not contained in $P$, 
this would imply that the slope $p/q$ of the annulus $E$ would satisfy $\vert p\vert=1$ if $n$ is even,
or that $\vert q\vert =1$ if $n$ is odd, contrary to the hypothesis.
So assume that $D$ and $E$ do not intersect. Let $F$ be a meridian disk of the solid torus $P$, this can 
intersect the arc $\alpha$ in many points, but suppose $\alpha$ has been isotoped so that its intersection 
with $F$ is minimal. The intersection of $F$ with $S^3- \Int H_8$ is a punctured disk, which we call $F$ again. 
Now look at the  intersections between $F$ and $D$. An innermost disk or outermost arc of intersection 
can be used to reduce  the number of points of intersection between $F$ and $\alpha$, which is not possible. 
So $D$ must be disjoint from $F$,  but this is not possible.
\end{proof}

Note that a knot $k$ which has a surface of type 8 may also have a surface of type 3.
Prolong the arc $\alpha$ on both ends until it touches $R_0$ in two points, but so that intersects
the $\smash{B_i^j}$'s only in two points, one lying in $E_1$, the other in $E_2$. 
Take curves $\gamma_1$ and $\gamma_2$ in $R_0$ of slope $p/q$, disjoint from the $C_i$'s, so that one
endpoint of $\alpha$ lies in $\gamma_1$ and the other in $\gamma_2$. 
$\Gamma_3 = \alpha\cup \gamma_1\cup \gamma_2$ is a $1$--complex as defined in \fullref{sec:type3}, and
$S'=\partial N(\Gamma_3)$ is disjoint from $S$, just take
$N(\gamma_3)$ thin enough. Note that if $n\geq 3$, the knot $k$ can be isotoped so that $k$ is contained in
$N(\Gamma_3)$; so if $k$ is well wrapped in $N(\Gamma_3)$, $S'$ will also be meridionally incompressible in
$S^3-k$. If $n=1,2$ then this construction may also work, depending if the knot $k$ can or cannot be isotoped
into
$N(\Gamma_3)$. However, if $n=2$ and the annuli $C_1$ and $C_2$ are nested, a similar construction yields a
surface of type 4, and if $n=1$, then we get a surface of type 1. In both cases the new surface may be 
meridionally incompressible.

Note that if $n$ is even, and the slope of the annulus $A_1$ is $1/q$, then with a little work it can be shown that 
the surface $S$ is in fact compressible. If $n$ is even, and the slope of $A_1$ is $p/1$, then the surface $S$ is 
isotopic to a surface with $n=2$. If $n$ is odd, and the slope of the annulus $A_1$ is $1/q$, then the surface
$S$ is isotopic to a surface with $n=2$. If $n$ is odd, and the slope of $A_1$ is $p/1$, then the surface $S$ is 
compressible.

Also note that if $n=2$, and the slope of the annulus $A_1$ is $p/1$, then both versions of a surface of type 8 are identical, ie the annuli $C_1$ and $C_2$ can be isotoped to be nested or non-nested.

\subsection{Surfaces of type 9}\label{sec:type9}

Let $R$ be a torus in $S^3$ constructed as follows. Let $A_1$, $A_2$, $\dots$, $A_n$ be $n$ annuli properly 
embedded in  $R_1$, all with slope $p/q$,  $\vert p\vert \geq 2$, $\vert q\vert \geq 2$ (in the usual 
coordinates of the solid torus
$R_0\cup (T\times [0,1])$). Suppose the annuli are nested, and say, $A_1$ is the innermost one. Let 
$\smash{B_1^1}$, $\smash{B_1^2}$, $\dots$, $\smash{B_n^1}$, $\smash{B_n^2}$ be $2n$ vertical annuli in $T\times I$, so that 
$\partial A_i \subset \partial (\smash{B_i^1\cup B_i^2})$, $1\leq i\leq n$. Let $C_1$, $\dots$, $C_n$ be $n$ 
annuli properly embedded in $R_0$, which are nested, whose boundaries coincide with the boundaries of 
the $\smash{B_i^j}$'s, so that $C_1$ is the innermost annulus and $\partial C_1\subset
\partial (\smash{B_1^2\cup B_2^2})$. Note that the union of the $A_i$'s, the $\smash{B_i^j}$'s and the $C_i$'s is an 
embedded torus, denoted by $R$. In the special case when $n=1$, assume that the annulus $C_1$ was chosen so that 
the torus $R$ does not bound a solid torus contained in $T\times I$. Note that in any case, $R$ is a 
standard torus in $S^3$.

Let $E$ be the annulus in $T_0$, with $\partial E\subset \partial (\smash{B_1^1\cup B_1^2})$, and whose interior is 
disjoint from the $\smash{B_i^j}$'s. Let $\gamma$ be a simple closed curve which is a core of the annulus $E$, and let
$N(\gamma)$ a small neighborhood of $\gamma$, disjoint from the $\smash{B_i^j}$'s. Let $\alpha$ be an arc in $T\times I$,
with one endpoint in
$\smash{B_1^1}$, the other in $N(\gamma)$, with interior disjoint from the $\smash{B_i^j}$'s, and so that it has a single local 
maximum in $T\times I$. By
sliding the arc $\alpha$, we can assume that its endpoints lie in the same level, say at level $\epsilon$, 
$0< \epsilon < 1$. We assume that the arc $\alpha$ cannot be isotoped, keeping its endpoints fixed, 
to a level arc lying in $T\times \{\epsilon\}$. 

Let $R'$ be the solid torus bounded by $R$ which does not contain the arc $\alpha$.
Let $H_9=R'\cup N(\alpha)\cup N(\gamma)$. This is a genus 2 handlebody; it can be seen as the solid tori 
$R'$ and $N(\gamma)$ joined by the $1$--handle
$N(\alpha)$, where $R'\cap N(\alpha)=B_1^1\cap N(\alpha)$ consists of a disk, say $E_1$, and 
$N(\alpha)\cap N(\gamma)$ consists of a disk $E_2$, which lies in level 
torus $T\times \{\epsilon\}$.
Let $S=\partial H_9$. We say that $S$ is a surface of type 9. 
See \fullref{fig:type8-9}, right.

Let $k$ be a knot in $H_9$, intersecting in two points each of 
$E_1$, $E_2$, and so that $k\cap N(\alpha)$ consists
of two arcs, each with a single local maximum, and $k\cap R'$ consists of one arc with a single local minimum, 
which is well wrapped in the solid torus $R'$. Furthermore, $k\cap N(\gamma)$ consists of one arc with a single
local minimum and which is well wrapped in $N(\gamma)$. 

\begin{thm}\label{thm:type9}
Let $S$ and $k$ be as above. $S$ is meridionally incompressible 
in $S^3-k$, 
and $k$ is a $(1,2)$--knot. 
\end{thm}

\begin{proof} By construction $k$ is in a $(1,2)$--position. We have to show that $S$ is incompressible in 
$S^3- \Int H_9$ and that $S$ is meridionally incompressible in $H_9-k$. It is not difficult to prove that $S$ is
meridionally incompressible in $H_9-k$. To do that consider the disks $E_1$, $E_2$, and look at the intersections 
between a compression disk $D$ and $E_1\cup E_2$. Using the hypothesis on $k$, we conclude that $D$ and 
$E_1\cup E_2$ can be made disjoint, which then implies that $D$ cannot exist. 

Suppose now that $S$ is compressible in $S^3- \Int H_9$, and let $D$ be a compression disk. 
Let $E'=E-\Int N(\gamma)$, these are two annuli. Look at the intersections between $D$ and the annuli $E'$. 
An argument as in the proof of \fullref{thm:type7} or \fullref{thm:type8} yields a contradiction.
\end{proof}

Note that a knot $k$ whose complement has a surface of type 9 may also have a surface of type 3.
The construction is identical to the one done in the previous section, just after the proof of
\fullref{thm:type8}.

Note that if $n$ is even, and the slope of the annulus $A_1$ is $1/q$, then the surface $S$ is isotopic to
a surface of type 7. If $n$ is even, and the slope of $A_1$ is $p/1$, then the surface $S$ is isotopic to a 
surface of type 3. If $n$ is odd, and the slope of the annulus $A_1$ is $1/q$, then the surface
$S$ is isotopic to a surface of type 3. If $n$ is odd, and the slope of $A_1$ is $p/1$, then the surface $S$ is 
isotopic to a surface of type 7.

\subsection[Incompressible tori and (1,2)-knots]{Incompressible tori and $(1,2)$--knots}\label{sec:tori}

It is not difficult to construct $(1,2)$--knots whose complement contains a closed meridionally incompressible
torus. There are three cases.

\begin{figure}[ht!]
\labellist
\small\hair 2pt
\pinlabel {$4$--braid} at 65 33
\pinlabel {$4$--braid} at 65 138
\pinlabel {$4$--braid in $A\times I$} at 223 130
\endlabellist
\centering
\includegraphics{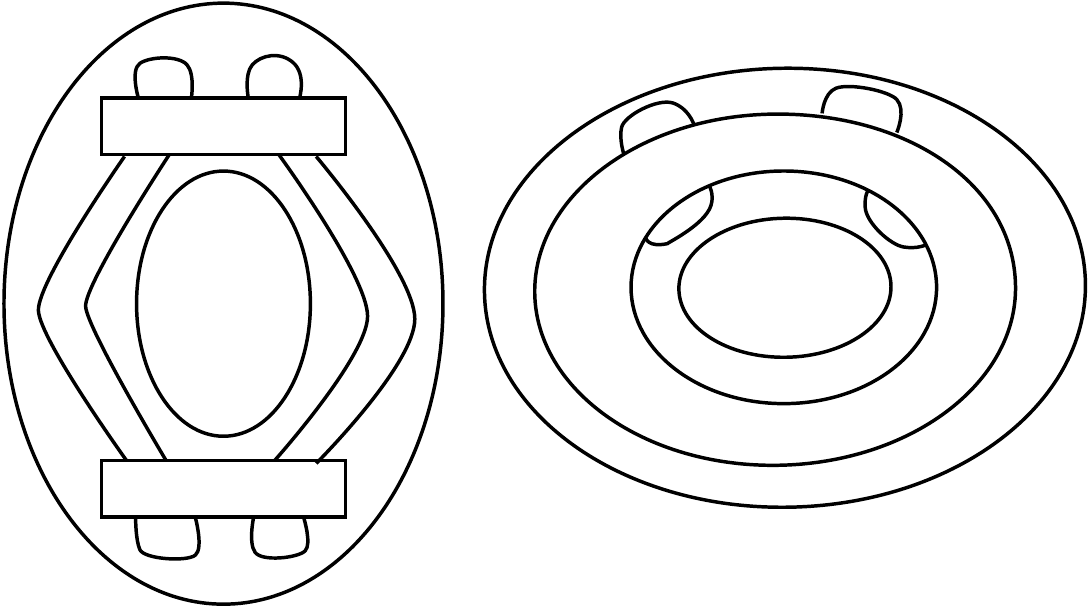}
\caption{}\label{fig:tori}
\end{figure}

Let $K$ be a nontrivial knot in a $(1,1)$--position. Embed a knot $k$ in $N(K)$ so that the wrapping 
number of $k$ in $N(K)$ is 2, and that $k$ is in a $(1,2)$--position. Clearly $S=\partial N(K)$ is
meridionally incompressible and $k$ is a $(1,2)$--knot. These knots inside $N(K)$ look like in
\fullref{fig:tori}, left.

Let $K$ be a nontrivial knot contained in the standard torus $T$ in $S^3$. Assume that $N(K)$ is of the
form $N(K)=A\times I \subset T\times I$, where $A$ is an annulus in $T$. Embed a knot $k$ in
$N(K)$ so that the wrapping  number of $k$ in $N(K)$ is $\geq 2$, and that $k$ is in a $(1,2)$--position. 
Clearly $S=\partial N(K)$ is meridionally incompressible and $k$ is a $(1,2)$--knot. 
These knots inside $N(K)$ look like in \fullref{fig:tori}, right.

Note that some of the knots constructed above could contain both, a meridionally incompressible 
surface of genus 2, and one of genus 1. Namely, let $S$ be a surface of type 1, and let $H$ be the
handlebody bounded by $S$. Then we can find a knot $K\subset H$, so that $K$ is in a $(1,1)$--position, 
and $S$ is incompressible in $H-K$ (but it is meridionally compressible). Now let $k$ be any $2$--cable 
of $K$ contained in $N(K)$. Then both surfaces $\partial N(K)$ and $S$ are meridionally incompressible 
in the complement of $k$. This construction may not work for the remaining types of surfaces, for it  
seems that there is no $(1,1)$--knot embedded in the regions bounded by that surfaces, so that the surfaces are
incompressible in the complement of that knot.

Let $S$ be a surface of type 7, and let $N$ be a big solid torus containing $S$, that is, a big 
neighborhood of the torus knot used in the construction of the surface of type 7. If the knot $k$
constructed in \fullref{sec:type7} is chosen so that $k$ lies in $N$, then both surfaces $S$ and 
$\partial N$ will be meridionally incompressible.

This construction, that is, of a big torus containing the surface, can also be done for surfaces of type 1, 2, 3, 4, and 6, if the surfaces and the knots are chosen adequately. But it may fail in the 
remaining cases. This is because it seems that a surface of type 5, 8 or 9 cannot be confined inside a solid torus
$N$, so  that $\partial N$ remains incompressible in $S^3-k$; except in the special case in surfaces of type 8, 
where $n=2$, and the lower annuli are non-nested. 

We showed before that knots with a surface of type 6 or 7 usually contain a surface of type 4, and knots
with a surface of type 8 usually contain a surface of type 3.
It follows that there are $(1,2)$--knots which contain 3 meridionally incompressible surfaces, 
one of genus 1, and two of genus 2 (one of type 4 and one of type 6 or 7, or one of type 3 and one of type 8).

\subsection{Further remarks}

All the knots $k$ constructed in this section have a $(1,2)$--presentation, that is, $b_1(k)\leq 2$.
So we could ask if they really have $b_1(k)=2$.  In \cite{E4}, all $(1,1)$--knots containing
a closed meridionally incompressible surface are described, and it is shown that the surfaces
are the boundary of a regular neighborhood of what is called a toroidal graph. So to show
that the present knots are not $(1,1)$--knots it suffices to show that the surfaces constructed here
do not satisfy the conditions given in \cite{E4}. This is clear for some of the surfaces of type 6 and 7,
the ones that do not bound a handlebody. It is intuitively obvious for the remaining cases, 
but a little more work is required to show that.

\section{Characterization of meridionally incompressible surfaces} \label{sec:characterization}

In this section we prove the following theorem.

\begin{thm}\label{thm:maintheorem} Let $K$ be a $(1,2)$--knot and let $S$ be a genus 2 meridionally 
incompressible surface in the complement of $K$. Then $K$ and $S$ come
from the construction of \fullref{sec:construction}, that is, $K$ and $S$ can be isotoped so that
$S$ looks as one of the surfaces of types 1, 2, 3, 4, 5, 6, 7, 8 and 9 constructed in
\fullref{sec:construction}.
\end{thm}

Let $T$ be a standard torus in $S^3$, and let $I=[0,1]$. Consider 
$T\times I\subset S^3$.
$T_0=T\times \{ 0\}$ bounds a solid torus $R_0$, and $T_1=T\times \{1\}$ 
bounds a solid torus $R_1$, such that $S^3= R_0\cup (T\times I)\cup R_1$. Let 
$k$ be a $(1,2)$--knot,  and assume that $k$ lies in $T\times I$, such that 
$k\cap T_0= k_0$ consists of two arcs,
$k\cap T_1 = k_1$ consists of two arcs, and $k\cap (T\times (0,1))$ 
consists of four vertical arcs.

Suppose there is a closed surface $S$ in $S^3-\Int N(k)$, which is 
incompressible and meridionally  incompressible.  Assume that
$S$ intersects $T_0$ and
$T_1$ transversely. Let $S_0=S\cap R_0$, $S_1=S\cap R_1$, and 
$\tilde S=S\cap (T\times I)$.
Let $\pi\co T\times I \rightarrow I$ be the height function, where we 
choose $0$ to be
the lowest point, and $1$ the highest. We may assume that the height 
function on
$\tilde S$ is a Morse function. So there is a finite set of different points
$X=\{x_1,x_2,\dots,x_m\}$ in $I$, so that
$\tilde S$ is tangent to $T\times \{ x_i \}$ at exactly one point, and this
singularity can be a local maximum, a local minimum, or a simple saddle. 
Suppose that $1 > x_1 > x_2 > \dots > x_m > 0$, 
that is, we numerate the singular
points starting from the upper level.
For any $y\notin X$,
$T\times \{ y \}$ intersects $\tilde S$ transversely, so for any such $y$,
$\tilde S \cap (T\times \{ y \})$ consists of a finite collection of simple closed
curves called level curves, and at a saddle point $x_i$, either one level 
curve of
$\smash{\tilde S}$ splits into two level curves, or two level curves are fused 
into one curve.

For example, any of the knots and surfaces constructed in 
\fullref{sec:construction} can be put in this position, after doing
an appropriate isotopy.

Define the complexity of $S$ by the pair 
$c(S)=(\vert S_0\vert + \vert S_1 \vert + \vert \tilde S \vert , \vert X \vert)$ 
(where $\vert Y \vert$ denotes the number of points if $Y$ is a finite set, 
or the number of connected components if it is a surface, and give to such 
pairs the lexicographical order). Assume that $S$ has been isotoped so that
$c(S)$ is minimal.

\begin{cla}\label{cla:claim1} The surfaces $S_0$, $S_1$ and $\tilde S$ are
incompressible and meridionally incompressible in $R_0$, $R_1$, and 
$T\times I - \Int N(k)$ respectively.
\end{cla}

\begin{proof} If  there is a meridian compression disk for one of the surfaces, 
then it will be also a meridian compression disk for $S$. 
Suppose then one of the surfaces is compressible, say $\tilde S$, and let
$D$ be a compression disk, which is disjoint from $k$. Then $\partial D$ is
essential in $\tilde S$ but inessential in $S$. By cutting $S$ along $D$ we 
get a surface $S'$ and a sphere $E$. Note that $S$ and $S'$ are isotopic in 
$M-k$.
For $S'$ we can similarly define the surfaces
$S_0'$, $S_1'$ and $\tilde S'$. Note that  
$\vert S_i \vert = \vert S_i' \vert + \vert E \cap R_i
\vert$, $i=1,\ 0$, then either $\vert S_0' \vert < \vert S_0 \vert$ or 
$\vert S_1' \vert < \vert S_1 \vert$, for $E$ intersects at least one of  $R_0$,
$R_1$. Also $\vert \tilde S' \vert \leq \vert \tilde S \vert$, so 
$c(S') < c(S)$,
but this contradicts the minimality of $c(S)$.
\end{proof}

This implies that $S_0$ is a collection of trivial disks, meridian
disks and incompressible annuli in $R_0$. If a component of $S_0$ is a trivial 
disk $E$, then $\partial E$ bounds a disk on $T_0$ which contains at least a 
component of $k_0$, for  otherwise $\vert S_0 \vert$ could be reduced. If a
component of $S_0$ is an incompressible annulus $A$, then $A$ is parallel to an
annulus $A^\prime \subset T_0$, and
$A^\prime$ must contain a component  of $k_0$, for otherwise $\vert S_0 \vert$ 
could be 
reduced.
Note that the slope of $\partial A$ can consist of one longitude and several meridians of $R_0$; 
in this case $A$ would also be parallel to $T_0-A^\prime$, and then the other 
component of $k_0$ would be in $T_0-A^\prime$. Note also that $S_0$ cannot
contain both incompressible annuli and meridian disks. A similar thing can be said for $S_1$.

\begin{cla}\label{cla:claim2} $\tilde S$ does not have any local maximum or minimum.
\end{cla}

\begin{proof} The proof is similar to that of 
\cite[Claim 2]{E4}. 
It consists in taking the maximum at the lowest level and then in pushing it down, getting
that either the surface is compressible or that $\tilde S$ has a component which is 
parallel to a subsurface in $T_0$. By pushing it into $R_0$ the complexity of $S$ is reduced.
\end{proof}

The proof of the claim also implies that if $S$ is in a position where $\tilde S$ has 
a maximum or a minimum, then $S$ can be isotoped to a position of lower
complexity.

Note that if at a certain nonsingular level $\{y\}$, there is a curve of intersection $\gamma$ 
which is trivial in the level torus $T\times \{y\}$, then $\gamma$ bounds a disk in the level torus which
intersects $k$ in two or more points, for otherwise $\tilde S$ will be compressible, meridionally 
compressible, or it would have a local maximum or minimum.

\begin{cla}\label{cla:claim3} Only the following types of saddle points are possible.
\begin{enumerate}
\item A saddle changing a trivial simple closed curve into two non-nested
trivial simple closed curves.
\item A saddle changing a trivial simple closed curve into two nested
trivial simple closed curves.
\item A saddle changing two non-nested trivial simple closed curves into a 
trivial simple closed curve.
\item A saddle changing two nested trivial simple closed curves into a 
trivial simple closed curve.
\item A saddle changing a trivial simple closed curve into two essential simple
closed curves.
\item A saddle changing two parallel essential curves into a trivial curve.
\item A saddle changing an essential curve $\gamma$ into a curve with the 
same slope
as $\gamma$, and a trivial curve.
\item A saddle changing an essential curve $\gamma$ and a trivial curve into an
essential curve with the same slope as $\gamma$.
\end{enumerate}
\end{cla}

\begin{proof} See \fullref{fig:sing1} and \fullref{fig:sing2}. At a saddle, either one level curve of
$\tilde S$ splits into two level curves, or two level curves are joined into one level curve.
If a level curve is trivial in the corresponding level torus and it bounds a 
disk intersecting $k$ in four points and at a saddle
the curve joins with itself, then the result must be either two non-nested 
trivial curves each bounding a disk intersecting $k$ in two points,  or two
essential simple closed curves, for otherwise $\tilde S$ would be compressible or
meridionally compressible; in this case the singularity is of type 1 or 5. 

If a level curve is trivial in the corresponding level torus and it bounds a 
disk intersecting $k$ in less than four points and at a saddle the curve joins
with itself, then the result must be either two essential simple closed curves,
or two nested trivial curves; in the latter case, the original curve must bound a
disk intersecting $k$ in two points, and one of the new curves bound a disk
intersecting $k$ in four points, for otherwise  $\tilde S$ would be compressible
or meridionally compressible. So we have a singularity of type 2 or 5. 

\begin{figure}[ht!]
\labellist
\small\hair 2pt
\pinlabel $1$ at 7 204
\pinlabel $2$ at 7 142
\pinlabel $3$ at 7 82
\pinlabel $4$ at 7 22
\endlabellist
\centering
\includegraphics[scale=.9]{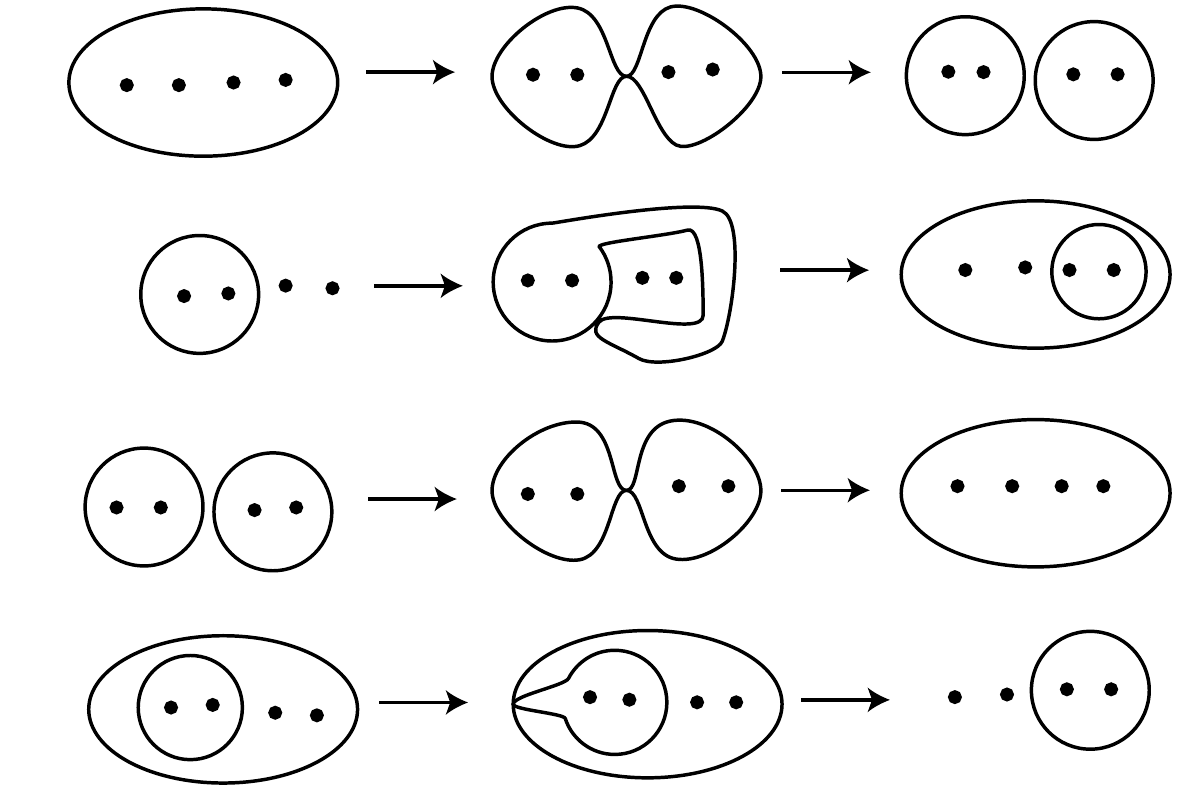}
\caption{}\label{fig:sing1}
\end{figure}

If two trivial level curves are joined into one and are non-nested, then each 
bounds a disk intersecting $k$ twice, and the new curve bounds a disk intersecting
$k$ in four points. This is a singularity of type 3. If two trivial level curves 
are joined into one and are nested,
then the innermost one bounds a disk intersecting $k$ twice, the outermost one
bounds a disk intersecting $k$ in four points, and the new one bounds a disk
intersecting $k$ twice. This gives a singularity of type 4.

If in a level there are nontrivial curves of intersection,
then there is an even number of them, for $S$ is separating. So if a curve is 
nontrivial and at the saddle joins with itself, then the result is a curve with the 
same slope as the original and a trivial curve, for the saddle must join points 
on the same side of the curve. This is a singularity of type 7.

Finally note that there may be singularities of types 6 and 8. 
\end{proof}

\begin{figure}[ht!]
\labellist
\small\hair 2pt
\pinlabel $5$ at 7 198
\pinlabel $6$ at 7 137
\pinlabel $7$ at 7 84
\pinlabel $8$ at 7 37
\endlabellist
\centering
\includegraphics[scale=.9]{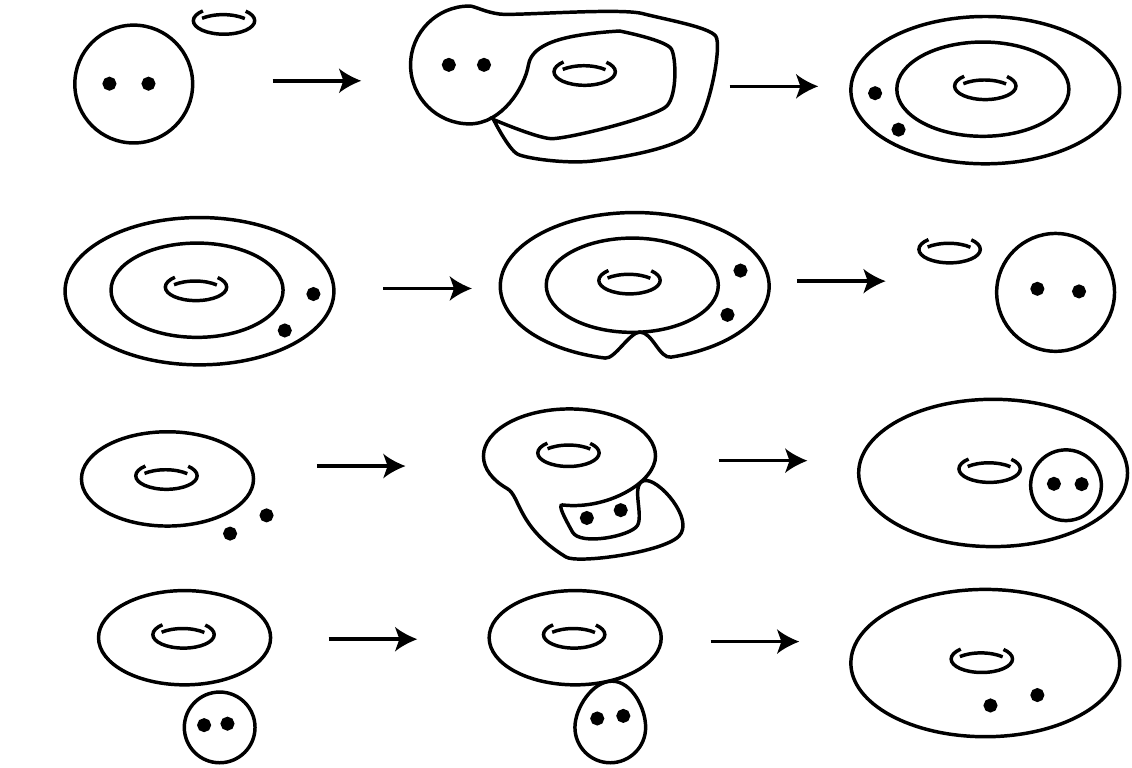}
\caption{}\label{fig:sing2}
\end{figure}

\begin{figure}[ht!]
\labellist
\small\hair 2pt
\pinlabel $1$ at 58 122
\pinlabel {$2$, $5$} at 177 122
\pinlabel $7$ at 290 122
\pinlabel $3$ at 58 7
\pinlabel {$4$, $6$} at 177 7
\pinlabel $8$ at 290 7
\endlabellist
\centering
\includegraphics[scale=.9]{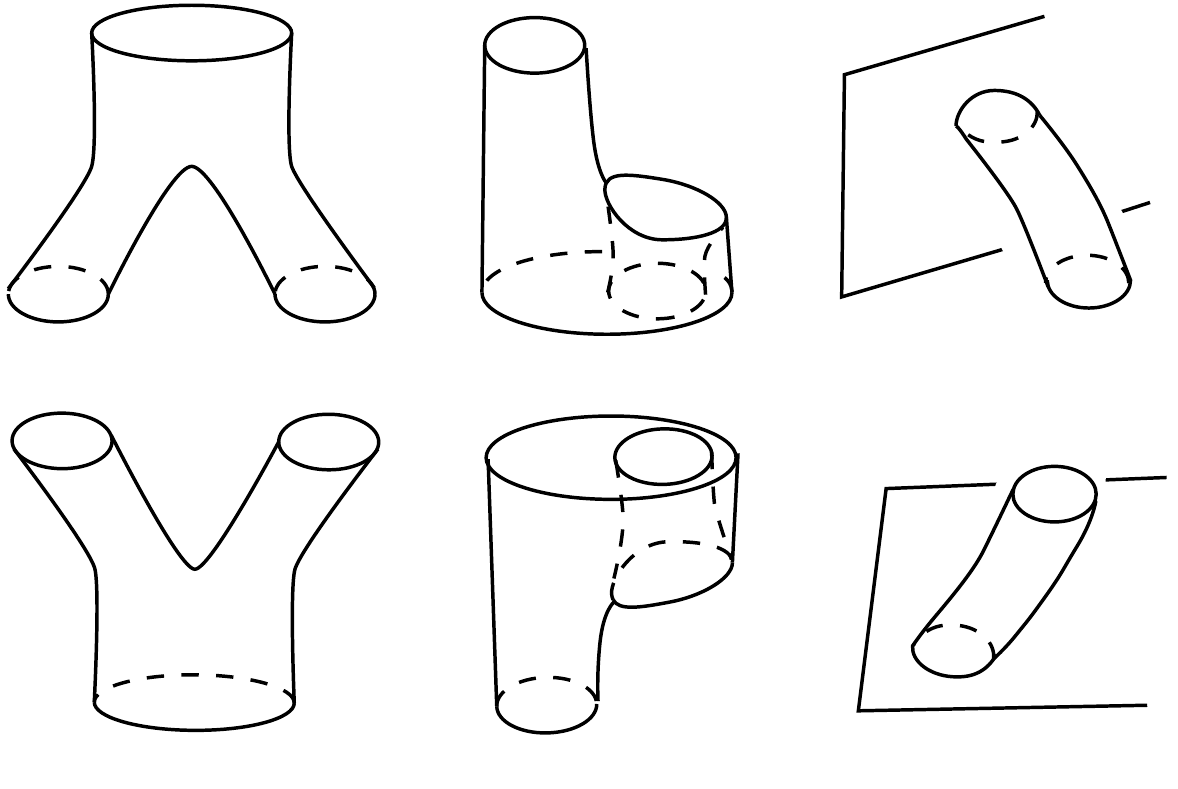}
\caption{}\label{fig:sing3}
\end{figure}

In \fullref{fig:sing3} we show locally how the surface $S$ looks in the neighborhood of a singularity. 
 
\begin{cla}\label{cla:claim4} Suppose that at a certain
non-singular level there is a curve $\gamma$ of intersection
which is trivial in the level torus and bounds a disk which
intersects $k$ in two points. If $\gamma$ is a trivial  curve
in $S$, then $\gamma$ bounds a disk $D \subset S$, which
consists of an annulus contained in $\tilde S$  with no singular
points and a disk component of $S_i$ which is trivial in $R_i$
($i=1,0)$.
\end{cla}

\begin{proof}  The curve $\gamma$ bounds a disk $E$ in the level
torus $T\times \{y\}$, and also bounds a disk $D\subset S$.
Suppose that $D$ is not as required. A collar neighborhood of
$\gamma$ in $D$ lies below $T\times\{y\}$, say. Suppose first
that $E\cap D=\partial E=\partial D$, but note that $E$ may
contain more curves of intersection with $S$. Now, $D\cup E$
bounds a $3$--ball $B$, which intersects $k$ in a spanning arc
$k'$. The $3$--ball $B$ can be isotoped to lie in a product
neighborhood of $E$, so that $k'$ is isotoped to an arc in this
product neighborhood, and it preserves its singular points.
Note that $D$ may contain disks of $S_0$ or $S_1$ which are
meridian disks of $R_0$ or $R_1$, and then an arc of $k$ may
wrap around a meridional annulus in $T_0$ or $T_1$; during the
isotopy, such arc must be arranged so that it now lies on a
trivial annulus contained in some level torus. Now pull $B$
down, eliminating singular points of $D$, until $D$ consists of
an annulus in $\tilde S$ without singular points and it
intersects $R_0$ in a trivial disk. If $int\, B \cap S$ is
non-empty, any component of $B\cap S$ must be a disk, for
otherwise a component of $\tilde S$, $S_0$ or $S_1$ would be
compressible. So any component of $B\cap S$ can be arranged to
be a disk consisting of an annulus in $\tilde S$  without
singular points and a trivial disk in $R_0$. At the end of this
procedure we have reduced $c(S)$. Now the arc $k'$ can be
rearranged so that the knot $k$ is in a $(1,2)$--position.

Suppose now that $E$ and $D$ intersect in some simple closed
curves. Look at the whole collection of curves $E\cap S$;
these curves are concentric in $E$, for each bounds a disk
intersecting $k$ twice, and each curve is trivial in $S$. Look
at the intersection curves in $S$, and among the innermost ones,
take the one which is outermost in $E$. Let $\alpha$
be this curve, which bounds a disk $D'\subset S$. As $D'\cap
E=\partial D'$, we can assume, by an argument as in the
previous paragraph, that $D'$ does not contain any singular
points of $\tilde S$ and that contain a single component of
$S_0$ or $S_1$. The curve $\alpha$ in $E$, with another curve of
intersection, say $\beta$  (perhaps $\beta=\gamma$), cobound an
annulus $E'$ in $E$, with interior disjoint from $S$. Consider
a copy of $D'$ and $E'$ to form a disk $D_1$ with interior
disjoint from $S$, and with $\partial D_1=\beta$. Cut $S$ with
$D_1$, getting a new surface $S'$ isotopic to $S$. If the curves
$\alpha$ and $\beta$ are concentric in $S$, the region bounded
by them in $S$ contains at least a component of $S_0$ or $S_1$,
which is eliminated in $S'$, so this reduces $c(S)$. If the
curves $\alpha$ and $\beta$ are non-concentric in $S$, the disk
$D''$ bounded by $\beta$ in $S$ must intersect the disk $E$ in
some simple closed curves and then must contain at least two
components of $S_0$ or $S_1$. By replacing $D''$ with $D_1$ we
reduce the complexity of $S$, for $D_1$ contains only a
component of $S_0$ or $S_1$. We may have introduced a new
local  maximum or minimum in $S'$ at level $\{y\}$. This new
singularity can be eliminated as in Claim $\ref{cla:claim2}$,
getting a surface with lower complexity.
\end{proof}

\begin{cla}\label{cla:claim5} Suppose that at a certain nonsingular level there is a curve $\gamma$ of 
intersection which is  trivial in the level torus and bounds a disk which intersects $k$ in four points, 
and its interior is disjoint from
$S$. If $\gamma$ is a trivial curve in $S$, then $\gamma$ bounds a disk $D \subset S$, which consists of an 
annulus contained in $\tilde S$ with no singular points and a disk component of $S_i$ trivial in $R_i$ ($i=1,0)$.
\end{cla}

\begin{proof} It is similar to the previous claim.
\end{proof}

\begin{cla}\label{cla:claim6} Suppose that at a certain level torus 
$T\times \{y \}$, there is a singularity $x_i$ of type 1, 2 or 7. 
So there is a curve $\gamma$ of intersection, which contains the singular point, it is trivial in the 
level torus, bounds a disk $E$, so that the other singular curve is not contained in $E$. Assume that 
$E$ intersects $k$ in two points. At a level just below
$y$, there is a curve $\gamma'$ of intersection, which is parallel to $\gamma$ in $S$, and which bound
a disk $E'$ in that level, also intersecting $k$ in two points. Then $\gamma'$ is a nontrivial curve in 
$S$.
\end{cla}

\begin{proof}
Suppose that $\gamma'$ is trivial in $S$. Then by \fullref{cla:claim4}, 
$\gamma'$, and in fact $\gamma$, bounds a disk 
$D \subset S$, which consists of an annulus contained in $\tilde S$ with no singular points and a disk 
component of $S_0$ trivial in $R_0$. Now, $E$ and $D$ bound a $3$--ball, which intersects
$k$ in a spanning arc $k'$. So $D$ can be isotoped into $E$, and then the singularity $x_i$ is eliminated. 
The arc $k'$ can be isotoped to lie in a product neighborhood of $E$, preserving its singular points, and now by
finding a vertical path from $T\times \{y\}$ to $T_0$ disjoint from $S$, $k'$ can be rearranged to be in a
$(1,2)$--position. This contradicts the minimality of $c(S)$.
\end{proof}

\begin{cla}\label{cla:claim7} Suppose that at a certain nonsingular level there is a curve $\gamma$ of 
intersection which is essential in the level torus. If $\gamma$ is a trivial curve in $S$, then $\gamma$  
is a meridian or a longitude of the level torus, and it bounds a disk $D \subset S$, which contains a 
meridian disk of $S_0$ or $S_1$.
\end{cla}

\begin{proof} Let $T\times \{y\}$ be the nonsingular level at which $\gamma$ lies.
Let $D$ be the disk in $S$ bounded by $\gamma$. The disk
$D$ intersects $T\times \{y\}$ in $\partial D$ and possibly in a collection of simple closed
curves, so the slope of $\partial D$ on $\gamma$, seen it as a knot in $S^3$, is the same as the slope of 
$(T\times \{y\})- \Int N(\gamma)$ in $\gamma$, so this has to be $0$, and then
$\gamma$ must be a meridian or a longitude
of $T\times \{y\}$. Now, $D$ must contain a meridian disk of $S_0$ or $S_1$, for otherwise it
will be contained in $S^3$ minus the cores of $R_0$ and $R_1$, but $\gamma$ is not trivial in that 
product region.
\end{proof}

\begin{cla}\label{cla:claim8} Suppose that at a certain nonsingular level, there are two 
concentric trivial curves of intersection, $\gamma_1$ and $\gamma_2$, which bound an annulus 
$A$ in such level, which is
disjoint from $k$. Suppose the curves are nontrivial in the surface $S$. Then the
curves cannot be parallel in  $S$.
\end{cla}

\begin{proof} Suppose the curves $\gamma_1$ and $\gamma_2$ are parallel in
$S$, then they bound an  annulus $A' \subset S$. Now $A$ and $A'$ bound a solid torus 
disjoint from $k$, and as
$\gamma_i$ is a trivial curve, it must be that $A$ and $A'$ are in fact isotopic.
If $A$ has more intersections with $S$, then we just have more annuli contained in
that solid torus, and could take an innermost one. The annulus $A'$ may have some singular points 
of $\tilde S$, if not then it contains an annulus of $S_0$ or $S_1$, which is parallel to an annulus in
$T_0$ or $T_1$ disjoint from $k$, which is not possible. So isotope $A'$ to $A$,
reducing the complexity of $S$, but putting it in a position, in which it has a maximum
or minimum. This can be eliminated as in \mbox{\fullref{cla:claim2}}, reducing then the complexity of
$S$.
\end{proof}

\begin{cla}\label{cla:claim9} Suppose that at a certain nonsingular level, there are two curves of 
intersection, $\gamma_1$ 
and $\gamma_2$, where $\gamma_1$ is trivial and $\gamma_2$ is essential in the level torus. 
If the curves are parallel in  $S$, then $\gamma_2$ is a meridian or longitude of the level torus, and
$\gamma_1\cup \gamma_2$ bounds an annulus containing a meridian disk of $R_0$ or $R_1$.
\end{cla}

\begin{proof}
If $\gamma_1$ and $\gamma_2$ are parallel in $S$, then there is an annulus
$A\subset S$, $\partial A =\gamma_1 \cup \gamma_2$. As $\gamma_1$ is trivial in the level torus, 
by an argument as in \fullref{cla:claim7} we see that this is possible 
only if $\gamma_2$ is a meridian or a longitude of the level torus, but then $A$ would intersect 
the core of $R_0$ or
$R_1$, implying that one of $S_0$ or $S_1$ must contain a meridian disk of $R_0$ or $R_1$.
\end{proof}

\begin{cla}\label{cla:claim10} Suppose that at a certain nonsingular level torus 
$T\times \{ y\}$, there are two curves of intersection, $\gamma_1$ and $\gamma_2$, which are essential in the 
level torus, and bound an annulus intersecting $k$ twice. Suppose that the 
curves are parallel in $S$, and cobound an annulus $A\subset S$ which is contained in $R_0\cup (T\times [0,y])$. 
Then either $A$ has no singular points of $S$ and contains just one annulus component of $S_0$, or it has  just a
type 6 singularity and contains just one disk component of $S_0$.
\end{cla}

\begin{proof} Suppose $A$ has singular points of $S$, and look at the first one. If it is of type 7, then 
it will contradict \fullref{cla:claim6}. So it must be of type 6, changing the two curves into a trivial
curve $\gamma_3$. Then by \fullref{cla:claim4} there cannot be more singularities in $A$, 
and $\gamma_3$ bounds a disk in $S$ which contains a disk of $S_0$. 
\end{proof}

\begin{cla}\label{cla:claim11} Suppose that at a certain nonsingular level torus 
$T\times \{ y\}$, there are two curves of intersection, $\gamma_1$ and $\gamma_2$, which are trivial 
in the level torus, are non-nested and each bounds a disk intersecting $k$ twice. Suppose that the 
curves are parallel in $S$, and cobound an annulus $A\subset S$ which is contained in $R_0\cup (T\times [0,y])$. 
Then $A$ has just a type 3 singularity and contains just one disk component of $S_0$.
\end{cla}

\begin{proof} Note that $A$ must have some singularities. The first one can be of type 2 or 3. 
If it is of type 2, it will contradict \fullref{cla:claim6}, so assume it is of type 3. So, after this
singularity we get a trivial curve $\gamma_3$ in a level torus, which bounds a disk intersecting $k$ in four
points, and $\gamma_3$ is trivial in $S$. By \fullref{cla:claim5}, 
$\gamma_3$ bounds a disk in $S$ which has
no singular points and contain one disk component of $S_0$, so the annulus is as desired.
\end{proof} 

Note that up to this point the arguments apply for any meridionally incompressible surface
in the complement of a $(1,2)$--knot. The next claims will make use of the hypothesis
that $\genus(S)=2$.

\begin{cla}\label{cla:claim12} Suppose that a certain nonsingular level there is a collection 
of concentric trivial curves, where the innermost and the outermost one bound an
annulus $A$ in such level, which is disjoint from $k$. Suppose the curves are
nontrivial in the surface $S$. Then the collection consists of at most 3 curves.
\end{cla}

\begin{proof} If there are more than 3 curves, then as they are nontrivial in 
$S$, and $S$ has genus 2, necessarily two of them will be parallel, contradicting
\fullref{cla:claim8}.
\end{proof}

These claims imply that most of the curves we see in a nonsingular level are essential curves in $S$.
 
\begin{cla}\label{cla:claim13} $S_1$ ($S_0$) does not contain both meridian disks
and trivial disks of $R_1$ ($R_0$).
\end{cla}

\begin{proof} Suppose $S_1$ contains meridian disks and trivial 
disks; it must consist of an even number of meridian disks, for it is
separating. If $\tilde S$ has no singular points, then it would be a collection of annuli,  
and $S$ would be a sphere; suppose then that there are some singularities on $\tilde S$.

If the first singularity is of type 2, then there is a trivial disk whose 
boundary joins with itself. Note that this singularity can be pushed to lie on $R_1$, 
changing the
trivial disk into an annulus bounded by trivial curves. One of these curves bounds a
disk $E$ intersecting $k$ in two points, corresponding to one of the arcs of $k$ lying
in $R_1$. As there are meridian disks, the arc can be rearranged to lie in $E$, but then $S$ 
will be compressible, as $\partial E$ must be essential in $S$ by 
\fullref{cla:claim6}.

If the first singularity is of type 3, then it can be pushed to lie 
on $R_1$, changing two trivial disks into a trivial disk, then reducing the complexity of $S$.

If the first singularity is of type 4, again it can be pushed to lie on $R_1$, changing two 
trivial disks into a trivial disk, then reducing the complexity of $S$.

If the first singularity is of type 5, then it can be pushed to lie on 
$R_1$, changing a trivial disk into an annulus $A$, whose boundary consists of curves parallel 
to the boundaries of the meridian disks. At least one of the arcs of $k_1$ must lie in the region 
between $A$ and $T_1$. So there is a disk $D$ in $R_1$, with $D\cap A=\partial D$, which is an 
essential curve on $A$. The remaining arc of $k_1$, say $k'$ can be arranged so that intersect $D$ 
in at most one point. If $k'$ intersects $D$ in one point then $S$ is meridionally compressible, 
and if it is disjoint from $D$ then $S$ is compressible, unless $\partial A$ is an inessential 
curve in $S$, but in this case by cutting $S$ with $D$ we get a surface $S'$ isotopic to $S$, 
but with $c(S') < c(S)$. 

If the first singularity is of type 6, then again it can be pushed to lie on $R_1$, changing two 
meridian disks into a trivial disk. The arcs of $k_1$ can be 
rearranged to be in the required position. So we have reduced the complexity of $S$.

If the first singularity is of type 7, then there is a meridian disk whose boundary joins with itself. 
Push this singularity to $R_1$, getting an annulus whose boundary consists of a meridian of $R_1$ and 
a trivial curve $\alpha_1$ which bounds a disk $E$ in $T_1$ intersecting $k_1$. Note that because there 
are trivial disks in $S_1$, the curve $\alpha_1$ bounds a disk $D$ in $R_1$ disjoint from $k_1$, and then $S$
will be compressible, unless $\alpha_1$ is an inessential curve in $S$, but this is not possible by
\fullref{cla:claim6}. 

If the first singularity is of type 8, then it can be pushed to lie on $R_1$, changing a meridian disk and 
a trivial disk into a meridian disk. Again, the arcs of $k_1$ can be rearranged to be in the required position.
So we have reduced the complexity of $S$.

Therefore assume that the first singularity is of type 1.
Then we have a collection of nested trivial disks, each bounding a disk intersecting $k$ in four points. 
Possibly we have a sequence of type 1 singularities, starting with the innermost curve. In each singularity 
the disk changes into an annulus. If there is a type 2, 3, 4, 5, 6, 7 or 8 singularity before all trivial disks
are transformed, then with the same arguments as above we see that there will be a compression disk for $S$, or
there is an isotopy which reduces the complexity of $S$. So we must have a sequence of type 1 singularities
which transforms the disks into a collection of nested annuli. Note that the new level curves obtained are
nontrivial in $S$, by \fullref{cla:claim6}. Then by \fullref{cla:claim12} there are at most 3 of these
singularities, so there are at most 3 nested annuli.  

Up to this point, the annuli can be seen as lying around an arc with one maximum in a region $A\times [y,1]$ 
(where $A \subset T_1$ is an annulus), together with a ball bounded
by two meridians disks of $R_1$, ie this region is just a $3$--ball $B'$. 
In the level torus $T\times \{y \}$ there are 3 sets of curves of intersection with $S$, consisting of curves 
$\alpha_i$'s, which are trivial and concentric in the level torus, essential and nonparallel in $S$; 
curves $\beta_i$'s, which are trivial and concentric in the level torus, essential and nonparallel in $S$, and
such that $\alpha_j$ is parallel to $\beta_j$ in $S$, for each $j$; curves $\gamma_i$'s, which are essential in
the level torus, are trivial in $S$, and bound a meridian disk in the solid torus 
$(T\times [y,1])\cup R_1$. Note also that the curves $\alpha_i$'s and 
$\beta_i$'s lie in the same component of $(T\times\{y\})-\cup\{\gamma_i\}$,
say in an annulus between $\gamma_1$ and $\gamma_2$.  

The next singularity could be of types 2, 3, 5, 6, 7 or 8. It is not difficult to see that if the next
singularity is of type 2, 3 or 5, then there will be a compression disk for $S$. If it is of type 6, then two
meridian disks are changed into a disk $E$. The path followed by the singularity can be complicated, but
because the nested annuli lie around an arc with just one maximum, these can be isotoped (in the $3$--ball
$B'$ determined by $A\times [y,1]$ and the $3$--ball bounded by two meridians disks  of $R_1$), so that the
path and the annuli look simple, and then it is not difficult to see that the disk $E$ can be isotoped so that
$E\cap R_1$ is a new trivial disk in $S_1$, reducing the complexity of $S$. If the next singularity is of
type 8, then it can be pushed into $R_1$, changing a meridian disk and a trivial disk into a meridian disk,
then reducing the complexity of $S$. 

So the only possibility left is that the next singularity is of type 7, where the path of the saddle encircles 
the nested annuli.
So a curve, say $\gamma_1$, is split into a nontrivial curve $\gamma_1'$ and a trivial curve $\beta_j$, which 
is concentric with the curves $\beta_i$'s, say. So if there are 3 nested annuli, we will have four parallel
curves, which contradicts \fullref{cla:claim12}. Suppose then that there are two nested annuli,
and then curves $\alpha_1,\alpha_2$, and $\beta_1,\beta_2$. After the singularity of type 7 we have one 
more curve, denoted $\beta_3$, and the curves $\beta_i$'s lie between $\gamma_1'$ and $\gamma_n$, say. 

So far, we have 3 singularities, $x_1$, $x_2$ of type 1 and $x_3$ of type 7. Look at the next singularity
$x_4$,  it may be of types 5, 6, 7 or 8. If it is of type 5, then the surface $S$ will be compressible.  If it
is of type 6 or 8, then such singularity can be pushed to a level above the singularity $x_3$.  So the
singularity 
$x_4$ is again of type 7. After this singularity we will have one more trivial curve of intersection,
which cannot be concentric with the $\beta_i$'s by \fullref{cla:claim12}, so it must be concentric with the
$\alpha_i$'s.  There are two cases, depending on which of $\gamma_1'$ or $\gamma_2$ is split by the
singularity. Suppose first  that $\gamma_1'$ is split into $\gamma_1''$ and $\alpha_3$. Note that
$\gamma_1''$ is the boundary of a disk 
$D$ (formed by a meridian disk bounded by $\gamma_2$ and an annulus between $\gamma_1''$ and $\gamma_2$), so $S$
is compressible, except if $\gamma_1''$ is trivial in $S$, but in that case, by cutting $S$ with $D$, we get a
new surface $S'$ isotopic to $S$, its embedding has a local minimum at $T\times[0,1]$, so that by isotoping it
we get a surface with lower complexity than $S$.

Suppose then that the singularity $x_4$ splits $\gamma_2$ into curves
$\gamma_2'$ and $\alpha_3$. If there are more than 2 curves $\gamma_i$, by the same arguments the next 
singularity would be of type 7, which again yields a contradiction with \fullref{cla:claim12}. 
So, up to this level, say
$y_1$, we have only the curves of intersection $\gamma_1',\gamma_2'$, $\alpha_1$, $\alpha_2$, $\alpha_3$,
$\beta_1$, $\beta_2$, $\beta_3$. Note that the curves $\gamma_1'$ and $\beta_3$ are parallel in $S$, and
$\gamma_2'$ and $\alpha_3$ are also parallel in $S$. Note that because $S$ is a surface of genus 2, these four
curves must be parallel in $S$ (for $\alpha_3$ ($\beta_3$) is nonparallel to $\alpha_1$ or $\alpha_2$
($\beta_1$ or $\beta_2$)). So in the solid torus 
$R_0\cup (T\times [0,y_1])$, there is an annulus $A$ bounded by two of these curves. 
In fact, $S\cap (R_0\cup (T \times [0,y_1]))$ consists of the annulus $A$ plus two pair of pants. Note that 
$\partial A$ cannot consist of $\gamma_i'$ and one of $\alpha_3$ or $\beta_3$, for $\gamma_i'$ is a longitude
of the solid torus 
$R_0\cup (T\times [0,y_1])$. If $\partial A=\alpha_3\cup\beta_3$, then as $A$ is a separating surface,
it leaves in one side the curves $\alpha_1$, $\alpha_2$, $\beta_1$, $\beta_2$, and in the  other side
the curves $\gamma_1'$, $\gamma_2'$, so there cannot be two pairs of pants with these six curves as their 
boundary. So we have $\partial A=\gamma_1' \cup \gamma_2'$, but this annulus is isotopic to an annulus in
$T\times \{y\}$,  by an isotopy that leaves fixed the knot and the rest of the surface. Then $c(S)$ was not
minimal.
 
Suppose then that $S_0$ contains just one trivial disk. 
The first singularity is of type 1, and after that level we have curves $\alpha_1$, $\beta_1$ and the 
$\gamma_i$'s. By the same arguments as above, the next two singularities must be of type 7. If there are 6 or
more curves $\gamma_i$, there will be 4 or more curves $\beta_i$, which is not possible. So suppose we have
just curves $\gamma_1$, $\gamma_2$, $\gamma_3$, and $\gamma_4$. By arguments as above, there are 4 singularities
of type 7, after which at a level $y_1$ we get curves $\{\alpha_1,\alpha_2,\alpha_3\}$,
$\{\beta_1,\beta_2,\beta_3\}$, and $\{\gamma_1',\gamma_2',\gamma_3',\gamma_4'\}$, where the following pair of
curves are parallel in $S$, $\{\alpha_1,\beta_1\}$, $\{\beta_2, \gamma_1'\}$, $\{\beta_3, \gamma_4'\}$,
$\{\alpha_2, \gamma_2'\}$, $\{\alpha_3,\gamma_3'\}$.  The curves $\alpha_i$'s and $\beta_i$'s lie in an annulus
in the level torus between the curves $\gamma_2'$ and $\gamma_3'$. Now in $S$ we have two sets of 4 parallel
curves, and $S\cap (R_0\times [0,y_1])$ consists of two annuli and two pair of pants. An argument as above
yields a contradiction.

Suppose now we have just two curves $\gamma_1$ and $\gamma_2$. After the singularity of type 1 we have 
curves $\alpha_1$, $\beta_1$, $\gamma_1$ and $\gamma_2$. Arguments as above show that the next two
singularities are of type 7. There are two possibilities. In one case, in a level $y_1$ just after
these singularities we get curves $\alpha_1$, $\beta_1$, $\beta_2$, $\beta_3$, $\gamma_1'$ and $\gamma_2'$,
where the pairs $\{\beta_2,\gamma_1'\}$ and $\{\beta_3,\gamma_2'\}$ are parallel curves in $S$. So in $S$ we
have 3 pairs of parallel curves, which implies that $S\cap (R_0\cup (T \times [0,y_1]))$ consists of two pairs of
pants, say $P_1$ and $P_2$. One of the pair of pants must have the curves $\gamma_1'$ and $\gamma_2'$. By the
position of the curves in the level torus $T\times \{y_1\}$, it is not possible that $\partial P_1=\gamma_1'\cup
\gamma_2'\cup \beta_1$ and $\partial P_2= \alpha_1\cup \beta_2\cup \beta_3$, so we must have $\partial
P_1=\gamma_1'\cup \gamma_2'\cup \alpha_1$ and $\partial P_2= \beta_1\cup \beta_2\cup \beta_3$. By the position
of the curves, we see that the curve $\beta_3$ must bound a disk disjoint from $S$ and disjoint from the knot,
so $S$ will be compressible.

The final case is that after the two type 7 singularities, at a level $y_1$, we have curves $\alpha_1$, 
$\alpha_2$, $\beta_1$, $\beta_2$, $\smash{\gamma_1'}$ and $\smash{\gamma_2'}$, where the $\alpha_i$'s and the $\beta_i$'s lie in
the same level annulus between $\smash{\gamma_1'}$ and $\smash{\gamma_2'}$, and the pairs of curves
$\{\alpha_2,\gamma_2'\}$, $\{\beta_2,\smash{\gamma_1'}\}$ are parallel in $S$.
If $\smash{\gamma_1'}$ and $\smash{\gamma_2'}$ are parallel in $S$, we will have a set of 4 parallel edges and an 
argument as above yields a contradiction. If they are nonparallel in $S$, then $S\cap (R_0\times [0,y_1])$
consists of two pairs of pants, say $P_1$ and $P_2$. One of the pair of pants must have the curves $\smash{\gamma_1'}$
and $\smash{\gamma_2'}$. So the only possibility (up to interchanging $\alpha_i$'s and $\beta_i$'s) is that $\partial
P_1=\smash{\gamma_1'\cup \gamma_2'}\cup \beta_1$ and $\partial P_2= \alpha_1\cup \alpha_2\cup \beta_2$, but by the
position of the curves this is not possible.
\end{proof}

\begin{cla}\label{cla:claim14} Suppose $S_1$ contains meridian disks of $R_1$, then $S_1$ contains 
exactly two meridian disks of $R_1$, and $S$ is a surface of type 5.
\end{cla}

\begin{proof} Suppose that $\partial S_1$ consists of curves $\gamma_1,\gamma_2,\dots,\gamma_n$, with $n$ an 
even number. Look at the first singularity, it may be of type 6 or 7, for there are no trivial disks in
$S_1$.  If it is of type 6, then two curves, say $\gamma_1$ and $\gamma_2$ are fused into a curve 
$\alpha_1$, which is trivial in a level torus and bounds a disk in such level torus intersecting $k$ in 2 or
4 points. But the curve $\gamma_1$ is trivial in $S$, so by 
\fullref{cla:claim4} or \fullref{cla:claim5},
$S$ can be isotoped to a position with lower complexity.

So assume the first singularity is of type 7. So a curve, say $\gamma_1$, is split into a nontrivial curve
$\gamma_1'$ and a trivial curve $\beta_1$, which bounds a disk in a level torus $T\times \{ y_1 \}$
intersecting $k$ in two points (for if it intersects $k$ in 3 or 4 points, there is a meridional
compression or compression disk for $S$, which is formed by the union of an annulus between $\gamma_1'$ and
$\gamma_2$ and a meridian disk bounded by $\gamma_2$). If the arcs of
$k_1$ lie in different components of
$T_1-S_1$, then again after passing the first singularity, we see that there is a compression disk for $S$. So
both arcs of $k_1$ lie in the annulus determined by $\gamma_1$ and $\gamma_2$, say. After the first
singularity,  the curve $\beta_1$ lies in an annulus determined by the curves $\gamma_1'$ and $\gamma_n$. The
second singularity may be of types 5, 6, 7 or 8. If it is of type 5, then the surface $S$ will be
compressible.  If it is of type 6 or 8, then the surface will be compressible or such singularity can be pushed
to a level above the first singularity. So the second singularity is again of type 7. 

Suppose first that the second singularity splits the curve $\gamma_1'$ into a nontrivial curve
$\gamma_1''$ and a trivial curve $\beta_2$. There are two cases, depending if the curves $\beta_1$ and
$\beta_2$ are concentric or not. Suppose first that the curves $\beta_1$ and $\beta_2$ are concentric.
They lie in a level annulus determined by $\gamma_1''$ and $\gamma_2$. Note that
$\gamma_1''$ is the boundary of a disk $D$ (formed by a meridian disk bounded by $\gamma_n$ and an annulus
between $\gamma_1''$ and $\gamma_n$), so $S$ is compressible, except if $\gamma_1''$ is trivial in $S$, but
in that case, by cutting $S$ with $D$, we get a new surface $S'$ isotopic to $S$, its embedding has a local
minimum at $T\times[0,1]$, so that by isotoping it we get a surface with lower complexity than $S$.

Suppose now that the curves $\beta_1$ and $\beta_2$ are nonconcentric.
They lie in a level annulus determined by $\gamma_1''$ and $\gamma_n$. 
Note that $\gamma_1''$ is the boundary of a disk $D$ (formed by a meridian disk bounded by $\gamma_2$ and an
annulus between $\gamma_1''$ and $\gamma_2$). So as before, $S$ is compressible or it can be isotoped to a
surface with lower complexity. 

Therefore the second singularity splits a curve other than $\gamma_1'$ into two curves. 
By these arguments, and as in the proof of the previous claim, we have a sequence of type 7 singularities
which create two sets of concentric trivial curves, $\alpha_i$'s and $\beta_i$'s. If $n\geq 8$,
there will in a certain level 4 of such concentric curves, contradicting \fullref{cla:claim12}. If $n=6$,
there will be in a  certain level torus $T\times \{y_1\}$, two sets of concentric trivial curves
$\{\alpha_1,\alpha_2,\alpha_3\}$,
$\{\beta_1,\beta_2,\beta_3\}$, plus 6 essential curves. Then in the surface $S$ we have two sets of 4 parallel 
edges. As in the proof of \fullref{cla:claim13}, this yields a contradiction. If $n=4$ there are two
possibilities, one is that in a level torus $T\times \{y_1\}$ there are two sets of concentric trivial curves
$\{\alpha_1,\alpha_2\}$,
$\{\beta_1,\beta_2\}$, plus 4 essential curves, or well two sets of concentric trivial curves $\{\alpha_1\}$,
$\{\beta_1,\beta_2,\beta_3\}$, plus 4 essential curves. 
An argument as in the proof of \fullref{cla:claim13} yields a contradiction.

So we must have $n=2$, and the first two singularities are of type 7. There are two possibilities for the curves 
obtained after these singularities, either we get two trivial concentric curves, or two trivial
curves which are nonconcentric, and in any case we get essential curves $\gamma_1'$ and $\gamma_2'$. Suppose
first that after the singularities we get trivial concentric curves
$\beta_1$  and $\beta_2$, plus the essential curves $\gamma_1'$ and $\gamma_2'$.
Note that $\beta_1$ and $\beta_2$ are essential and nonparallel in $S$, by \fullref{cla:claim6} and 
\fullref{cla:claim8}.
The third singularity can be of types 2, 5, 6, 7 or 8. 
Let $T\times \{y_1\}$ be a level torus just below this singularity,
and let $W=R_0 \cup (T\times [0,y_1])$, this is a solid torus.
If the third singularity is of type 2, the curve $\beta_2$ is split into
curves $\alpha_1$ and $\alpha_2$, where $\alpha_1$ bounds a disk intersecting $k$ in two points and
$\alpha_2$ encircles $\beta_1$ and $\alpha_1$; the curve $\alpha_1$ is essential in $S$ by
\mbox{\fullref{cla:claim6}}, and the curve $\alpha_2$ is also essential, for otherwise by \fullref{cla:claim4}, 
$\alpha_2$ bounds a disk $D$ with no
singularities, which then implies that $\gamma_1'$ and $\gamma_2'$ cobound an annulus which can be isotoped
into $T\times I$, contradicting the minimality of $c(S)$. The surface $S\cap W$ consists of an annulus and a
pair of pants. Because of the position of the curves, the annulus bounds the curves $\alpha_1$ and $\beta_1$,
and the pair of pants has as boundary the curves $\gamma_1'$, $\gamma_2'$ and $\alpha_2$. So, by \fullref{cla:claim11}, in the annulus there must be a singularity of type 3, but this singularity can be
interchanged with the type 2 singularity, to become a singularity between
$\beta_1$ and $\beta_2$, which implies that $S$ is compressible. 

If the third singularity is of type 5, the curve $\beta_2$ is
split into curves $\alpha_1$ and $\alpha_2$, which are parallel in the level torus $T\times \{y_1\}$ 
to $\gamma_1'$, 
and which are essential in $S$, for otherwise $\alpha_1$ or $\alpha_2$ will bound a disk contained in $W$,
which is not possible. The surface
$S\cap W$ consists of an annulus and a pair of pants. Because of the position of the curves, the
annulus bounds, say, the curves $\alpha_1$ and $\gamma_1'$, and the pair of pants has as boundary the
curves $\gamma_2'$, $\alpha_2$ and $\beta_1$. So in the pair of pants there must be a singularity of type 8,
joining $\beta_1$ and $\alpha_2$, but this singularity can be interchanged with the type 5 singularity, to
become a singularity between $\beta_1$ and $\beta_2$, which implies 
that $S$ is compressible. 

If the third singularity is of type 6, the curves $\gamma_1'$ and $\gamma_2'$ are fused into a curve $\alpha_1$, 
which is essential in $S$, and in the level torus $T\times \{y_1\}$ cobound an annulus with $\beta_2$ which 
intersects $k$ in two points. Now the curves $\alpha_1$, $\beta_1$ and $\beta_2$ must bound a pair of pants 
in the solid torus $W$. Then the curve $\alpha_1$ bounds a disk disjoint from $k$, so $S$ is compressible. 

If the third singularity is of type 7, the curve $\gamma_1'$ (or $\gamma_2'$) is split into a curve $\gamma_1''$
and a trivial curve $\alpha_1$. Again in $S\cap W$, there is a pair of pants and an annulus bounded by the
curves $\gamma_1''$, $\gamma_2'$, $\alpha_1$, $\beta_1$ and $\beta_2$. The only possibility is that the pair of
pants bounds $\gamma_1''$, $\gamma_2'$ and $\beta_2$, and the annulus bounds $\alpha_1$ and $\beta_1$, so
$\alpha_1$ is nonconcentric with $\beta_2$ in the level torus, by 
\fullref{cla:claim8}. Note that the pair of pants can be
isotoped into the level torus $T\times \{y_1\}$, so that $c(S)$ can be reduced.

If the third singularity is of type 8, the curve $\beta_2$ and $\gamma_1'$ (or $\gamma_2'$) are fused into 
a curve $\gamma_1''$. The surface $S\cap W$ consists of a pair of pants with boundary curves $\gamma_1''$,
$\gamma_2'$ and $\beta_1$. Then there must be one more singularity of type 8. If it is between 
$\beta_1$ and $\gamma_1''$, it can be interchanged with the previous singularity, getting then a singularity
between $\beta_1$ and $\beta_2$, which shows that $S$ is compressible. If the singularity is between
$\beta_1$ and $\gamma_2'$, then after isotoping $S$ in $R_0$, the singularities can be interchanged,
showing that $S$ is compressible.
 
So assume that after the first two type 7 singularities we have two nonconcentric trivial curves 
$\alpha_1$ and $\beta_1$, and two essential curves $\gamma_1'$ and $\gamma_2'$. Note that the pairs of
curves $\{\alpha_1,\gamma_2'\}$ and $\{\beta_1,\gamma_1'\}$ are parallel in $S$, and that curves from
different pairs are nonparallel in $S$; this is because if they are parallel, then there is an annulus 
cobounded by $\alpha_1$ and $\beta_1$, and a twice punctured torus bounded by $\gamma_1'$ and $\gamma_2'$, 
but in this case $S$ would be compressible. The next singularity can be of types 2, 3, 5, 6, 7 or 8. 
Let $T\times \{y_1\}$ be a level  torus just below this singularity, and let $W=R_0 \cup (T\times [0,y_1])$.

If the third singularity is of type 3, the curves $\alpha_1$ and $\beta_1$ are fused into a curve
$\alpha_2$ which bounds a disk in the level torus $T\times \{y_1\}$ intersecting $k$ in four points. 
One possibility is that the surface $S\cap W$ is a pair of pants bounded by $\gamma_1'$, 
$\gamma_2'$ and $\alpha_2$. Note that such a pair
of pants can be isotoped into $T\times \{y_1\}$, then $c(S)$ is not minimal. The other possibility is that
$S\cap W$ consists of a once-punctured torus bounded by $\alpha_2$ and an annulus bounded by 
$\gamma_1'$ and $\gamma_2'$. But note that because the curves $\gamma_1'$ and $\gamma_2'$ are
longitudinal in $W$, such an annulus bounded by $\gamma_1'$ and $\gamma_2'$ can be isotoped into 
$T\times \{y_1\}$, reducing $c(S)$.

If the third singularity is of type 7, the curve $\gamma_1'$ is split into $\gamma_1''$ and $\beta_2$,
where now $\beta_2$ is concentric with $\beta_1$. After that level $S\cap W$ consists of an annulus and 
a pair of pants, and by the configuration of the curves this is not possible.

If the third singularity is of type 6,
the curves $\gamma_1'$ and $\gamma_2'$ are fused into a curve $\alpha_2$ which is trivial in the level torus
and encircles $\alpha_1$ and $\beta_1$. Below that level, there are two possibilities for the surface $S\cap W$.
If it is a pair of pants bounded by $\alpha_1$, $\beta_1$ and $\alpha_2$, then it can be isotoped to the level
torus, and so $c(S)$ is not minimal. The other case is that $S\cap W$ is an annulus bounded by
$\alpha_2$ and $\alpha_1$ and a once-punctured torus bounded by $\beta_1$.
In that case, the fourth singularity is of type 4, splitting $\beta_1$ into $\beta_2$ and $\beta_3$. 
Note that the third and fourth singularities can be interchanged, and the type 6 singularity can be pushed to
$R_0$, reducing $c(S)$.

If the third singularity is of type 8, the curves $\gamma_1'$ and $\beta_1$, say, are fused into a curve
$\gamma_1''$. After that level the surface $S\cap W$ is a pair of pants or an annulus and a once-punctured
torus. If it is a pair of pants, then it can be isotoped to the level torus, contradicting the minimality of
$c(S)$. So, it consists of an annulus and a once-punctured torus. Then the once-punctured torus is bounded
by $\alpha_1$ and the annulus by $\gamma_1''$ and $\gamma_2'$. So there must be a singularity of type 5, 
and the part of $S$ below a certain level torus $T\times \{y_2\}$ consists of two annuli. By \fullref{cla:claim10}, $S_0$ consists of two annuli (which are isotopic to nested annuli, and anyone of them could
be the innermost one), or there are two more singularities of type 6 and $S_0$ consists of two nested trivial 
disks (again, anyone of them could be the innermost one).
If there is just one more singularity of type 6, then it can be pushed into $R_0$, reducing $c(S)$. Note that the
singularities of  type 8 and 5 can be interchanged, and that the singularity of type 8 and one of the singularities
of type 6 can be interchanged, giving one of type 6 and one of type 4. In any case we get a surface of type 5.

If the third singularity is of type 2, then similar arguments show that there is one more singularity of type
4, followed by one of type 8, concluding with a trivial disk and an annulus in $R_0$, or well, there is one more
singularity of type 6, and $S_0$ consists of two nested trivial disks. Again, the singularities of type 8 and 
6 can be interchanged, giving one of type 6 and one of type 4. In any case this is a surface of type 5.

If the third singularity is of type 5, then similar arguments show that there is one more singularity of type 8, 
concluding with two annuli or two disks in $R_0$. 
This is a surface of type 5.

Summarizing, the following sequences of singularity types are possible, all producing a surface of type 5.

(a)\qua 7, 7, 2, 4, 8, and $S_0$ consists of a trivial disk and an annulus, which are nested.

(b)\qua 7, 7, 2, 4, 8, 6, or 7, 7, 2, 4, 6, 4,  and $S_0$ consists of two nested trivial disks.

(c)\qua 7, 7, 8, 5, or 7, 7, 5, 8, and $S_0$ consists of two nested annuli.

(d)\qua 7, 7, 8, 5, 6, 6, or 7, 7, 5, 8, 6, 6, or  7, 7, 5, 6, 8, 6, or 7, 7, 5, 6, 6, 4, and $S_0$ consists 
of two nested trivial disks.
\end{proof}

Note that the proof of this claim also shows that if $S_1$ contains meridian disks, then $S_0$ has two 
components, two disks, two annuli or a disk and an annulus.

\begin{cla}\label{cla:claim15} Suppose that $S_0$ and $S_1$ consist only of trivial disks and annuli. 
Suppose that at a certain  nonsingular level there is a curve $\gamma$ of intersection which is essential 
in the level torus.  Then $\gamma$ is an essential curve in $S$.
\end{cla}

\begin{proof} This follows from \fullref{cla:claim7}.
\end{proof}

\begin{cla}\label{cla:claim16} Suppose $S_1$ consists of just one trivial disk. 
Then $S$ is a surface of type 1, 2, 3 or 4.
\end{cla}
 
\begin{proof} As $S_1$ contains a single disk, $S_0$ cannot contain meridian disks by the remark following \fullref{cla:claim14}. So by \fullref{cla:claim15}, any level curve of intersection which is essential in
a level torus, it must be essential in $S$.

Let $D$ be the disk component of $S_1$, then $\partial D$ bounds a disk in $T_1$ 
which contains both arcs of $k_1$, for otherwise $S$ will be nonseparating. If the
first singularity is of type 5, then it can be pushed to lie in $R_1$, changing $D$
into an annulus, but this reduces the complexity of $S$, which is not possible. So
the first singularity must be of type  1, which splits the trivial curve into two
trivial curves $\gamma_1$ and $\gamma_2$, each bounding a disk intersecting $k$ in two points. 

The next singularity has to be of type 3 or of type 2 or 5. Suppose first that it is of type 3.
Then the curves $\gamma_1$ and $\gamma_2$ are fused into a trivial curve $\gamma_3$, which bounds
a disk in a level torus intersecting $k$ in four points. Up to this level the surface obtained is a once-punctured torus. 
The next singularity can be of type 1 or 5. Suppose it is of type 5; so $\gamma_3$ is split into two curves 
$\gamma_4$ and $\gamma_5$, which are essential in the corresponding level torus. These curves are essential in $S$,
so they must be parallel in $S$, for otherwise the genus of $S$ is greater than 2. Then by \fullref{cla:claim10}, we can assume  that $\gamma_4$ and $\gamma_5$ lie on $T_0$ and bound an annulus
component of $S_0$. Then $S$ is as surface of type 1, as defined in 
\fullref{sec:type1}. 

Suppose now the third singularity is of type 1, so $\gamma_3$ is split into two trivial curves
$\gamma_4$ and $\gamma_5$, each bounding a disk in a level torus intersecting $k$ in two points.
Note that $\gamma_4$ and $\gamma_5$ are essential in $S$ by \mbox{\fullref{cla:claim6}}, and that $\gamma_4$ and
$\gamma_5$ must be parallel in $S$, because $\genus(S)=2$. Then by \fullref{cla:claim11}, the next singular
point must be of type 3, fusing 
$\gamma_4$ and $\gamma_5$
into a trivial curve $\gamma_6$. So, $\gamma_6$ must lie at $T_0$, and it bounds a disk in $R_0$.
So, $S$ is a surface of type 1.

Suppose now that the second singularity is of type 5. So, the curve $\gamma_1$ splits into two curves
$\gamma_3$ and $\gamma_4$, which are essential in the corresponding level torus.
The third singularity can be of type 5, 6, 7 or 8.

If the third singularity is of type 8, say fusing $\gamma_2$ and $\gamma_3$, then the second and third 
singularities can be interchanged, so that the new second singularity is of type 3 and we are in the previous case.

If the third singularity is of type 5, then $\gamma_2$ is split into two curves $\gamma_5$ and $\gamma_6$, 
which are essential in a level torus, say $T\times\{y_1\}$. Then as $\gamma_3$, $\gamma_4$, $\gamma_5$ and 
$\gamma_6$ are essential in $S$, among them there must be two pairs of parallel curves of $S$, which bound annuli
$A_1$ and $A_2$ in $S$. So $A_1$ and $A_2$ are contained in $R_0\cup (T \times [0,y_1])$. 
If $A_1$ bounds $\gamma_3$ and $\gamma_6$ and $A_2$ bounds $\gamma_4$ and $\gamma_5$, then one of the annuli,
say $A_1$, can be isotoped to lie in $T\times I$, reducing $c(S)$.
So assume that $A_1$ bounds $\gamma_3$ and $\gamma_4$ and that $A_2$ bounds
$\gamma_5$ and $\gamma_6$. By \fullref{cla:claim10}, we can assume
that either there are no more singularities in $S$ and that $S_0$ consists of two annuli, or that there
two more singularities of type 6 and $S_0$ consists of two disks (which have to be nested, for otherwise
$c(S)$ could be reduced). Note also that if there is just one more singularity of type 6, then it can be 
pushed into $R_0$, reducing $c(S)$. If the annuli $A_1$ and $A_2$ are non-nested in 
$R_0\cup (T\times [0,y_1])$ we have a surface of type 3, and if the annuli are nested then we have a 
surface of type 4. 

If the third singularity is of type 6. Then the two curves $\gamma_3$ and $\gamma_4$ are
fused into a trivial curve $\gamma_5$. The surface up to that level is a twice punctured torus.
The curve $\gamma_2$ is essential in $S$ by \fullref{cla:claim6}. 
If $\gamma_5$ is essential in $S$, then $\gamma_2$ and $\gamma_5$ must be parallel in $S$. By \fullref{cla:claim11}, the next singularity must
be of type 3, fusing $\gamma_2$ and $\gamma_5$ into a trivial curve $\gamma_6$ which must bound a disk in
$R_0$.  So we have a surface of type 2. If $\gamma_5$ is trivial in $S$, then it bounds a disk in $S$ with no
more  singularities by \fullref{cla:claim4}. It follows that the next singularity must be of type
5, changing $\gamma_2$ into two curves $\gamma_6$ and $\gamma_7$. These curves must be parallel in $S$,
so that we may assume that in a level torus $T\times\{y_1\}$ there are 3 curves, $\gamma_6$, $\gamma_7$ and 
$\gamma_5$. The part of $S$ contained in $R_0\cup (T\times [0,y_1])$ consists of a disk and an annulus.
If the disk and the annulus are non-nested in $R_0\cup (T\times [0,y_1])$, then we can assume that $S_0$
consists  of a trivial disk and an annulus and we have a surface of type 3. If the disk and annulus are 
nested then we  have a surface of type 4, but in this case there may be one more singularity of type 6, 
so $S_0$ consists  of a disk and one annulus, or of two disks.

Finally suppose that the third singularity is of type 7. Then the curve $\gamma_3$, say, is split into 
two curves $\gamma_5$ and $\gamma_6$, where $\gamma_5$ is concentric with $\gamma_2$ and $\gamma_6$
is parallel to $\gamma_4$. The curve $\gamma_5$ is nontrivial in $S$, by \fullref{cla:claim6}. Also, the
curves $\gamma_2$ and $\gamma_5$ are nonparallel in $S$ by \fullref{cla:claim8}. The curves $\gamma_4$ and
$\gamma_6$  are also nontrivial in $S$ by \fullref{cla:claim15}. The curve $\gamma_4$ (or $\gamma_6$) cannot
be parallel in $S$  to $\gamma_3$ or $\gamma_5$, by \fullref{cla:claim9}. Then $\gamma_4$ and $\gamma_6$
must be parallel in $S$, for otherwise $\genus(S) > 2$. Note that the surface up to a level just
below the third singularity, union with the annulus bounded by $\gamma_4$ and $\gamma_6$ is a twice punctured
torus, so that for $S$ to be a genus 2 surface, the curves $\gamma_2$ and $\gamma_5$ must be parallel in $S$,
a contradiction.

Suppose now that the second singularity is of type 2. So, the curve $\gamma_1$ splits into two curves
$\gamma_3$ and $\gamma_4$, which are trivial and concentric with $\gamma_2$ in the corresponding level torus.
The third singularity can be of type 4 or 5. If it is of type 4, then the curves $\gamma_3$ and $\gamma_4$
are fused into a curve $\gamma_5$ which is trivial in the corresponding level torus.
The surface up to that level is a twice punctured torus.
The curve $\gamma_2$ is essential in $S$. If $\gamma_5$ is essential in $S$,
then $\gamma_2$ and $\gamma_5$ must be parallel in $S$. So the next singularity must be of type 3,
fusing $\gamma_2$ and $\gamma_5$ into a trivial curve $\gamma_6$ which must bound a disk in $R_0$. 
So we have a surface of type 2. If $\gamma_5$ is trivial in $S$, then the next singularity must be of type
5, changing $\gamma_2$ into two curves $\gamma_6$ and $\gamma_7$. These curves must be parallel in $S$.
So in a level torus $T\times \{y_1\}$ there are 3 curves, $\gamma_6$, $\gamma_7$ and $\gamma_5$, and the part of 
$S$ below that level consists of a trivial disk and an annulus. If the annulus and the disk are non-nested we
have a surface of type 3, and if they are nested, we have a surface of type 4, but in this last case we may 
have one more singularity of type 6.  Finally, suppose the third singularity is of type 5. In this case the curve
$\gamma_4$ transforms into two curves $\gamma_5$ and $\gamma_6$, which are essential in the corresponding level
torus. This situation is identical to the situation in the preceding paragraph, so it is not possible.

Summarizing, the following sequences of singularity types are possible:

(a)\qua 1, 3, 5, or 1, 5, 8, $S_0$ consists of an annulus, and $S$ is a surface of type 1.

(b)\qua 1, 3, 1, 3, $S_0$ consists of a trivial disk, and $S$ is a surface of type 1.

(c)\qua 1, 5, 5, $S_0$ consists of two non-nested annuli, and $S$ is a surface of type 3.

(d)\qua 1, 5, 5, or 1, 5, 5, 6, 6, $S_0$ consists of two nested annuli, or of two nested disks, 
and $S$ is a surface of type 4.

(e)\qua 1, 5, 6, 3, $S_0$ consists of a trivial disk, and $S$ is a surface of type 2.

(f)\qua 1, 5, 6, 5, $S_0$ consists of a trivial disk and an annulus, which are non-nested, 
and $S$ is a surface of type 3.

(g)\qua 1, 5, 6, 5, or 1, 5, 6, 5, 6, $S_0$ consists of a trivial disk and an annulus, which are nested, or of
two nested trivial disks, and $S$ is a surface of type 4.

(h)\qua 1, 2, 4, 3, $S_0$ consists of a trivial disk, and $S$ is a surface of type 2.

(i)\qua 1, 2, 4, 5, $S_0$ consists of a trivial disk and an annulus, which are non-nested and 
$S$ is a surface of type 3.

(j)\qua 1, 2, 4, 5, or 1, 2, 4, 5, 6, $S_0$ consists of a trivial disk and an annulus, which are nested, 
or of two nested trivial disks, and $S$ is a surface of type 4.
\end{proof}

\begin{cla}\label{cla:claim17} Suppose $S_1$ consists of two trivial disks. Then $S$ is a surface of type 4, 5, 6 or 7.
\end{cla}

\begin{proof}
Suppose that the disk components of $S_1$ are non-nested, so each curve bounds a disk in $T_1$ containing one 
arc of $k_1$. The first singularity must be of type 2, 3 or 5. Note that in any case
such singularity can be pushed into $R_1$, showing that
$S$ is compressible, or well, that a disk is changed by an annulus, or that the two disks are 
changed by one disk. 
In any case the complexity of $S$ is reduced. So suppose the disks are nested. 

Suppose first that each of these disks bounds a disk in $T_1$
containing the two arcs of $k_1$. So just below $T_1$, we have two concentric curves bounding a disk
which intersects $k$ in 4 points. The first two singularities must be of type 1, so just after
these singularities we have two pairs of parallel curves, say $\{\gamma_1, \gamma_2\}$ and $\{\gamma_3,\gamma_4\}$,
each bounding a disk intersecting
$k$ in two points. So up to this level the surface is just two nested annuli. The third singularity must 
be of type 2, 3 or 5. Suppose first it is of type 5. So, say, the curve $\gamma_1$ is split into two
curves $\gamma_5$ and $\gamma_6$, which are essential in the corresponding level torus. Note that the curves
$\gamma_2$, $\gamma_5$ and $\gamma_6$ are nontrivial  in $S$ (by \fullref{cla:claim6} and \fullref{cla:claim15}). 
By \fullref{cla:claim8} the curves $\gamma_3$ and $\gamma_4$ are 
nonparallel in $S$, and note also that the
curves $\gamma_5$ and $\gamma_6$ cannot be parallel in $S$ to one of the curves $\gamma_3$ or
$\gamma_4$, by \fullref{cla:claim9}. So $\gamma_5$ and $\gamma_6$ must be parallel in $S$ for otherwise the
genus of $S$  is $> 2$. So $\gamma_5$ and $\gamma_6$ cobound an annulus in $S$. Then there should be a 
singularity of type 5 joining the curve $\gamma_2$ with itself. This shows that the genus  of $S$ is 
$> 2$.

Suppose now that the third singularity is of type 2, so, say, the curve $\gamma_1$ is split into two
curves $\gamma_5$ and $\gamma_6$, which are trivial in the corresponding level torus. 
The next singularity must be of type 5, changing the curve $\gamma_6$, say, into curves $\gamma_7$ and
$\gamma_8$, which are essential in a level torus. Note that $\gamma_3$, $\gamma_4$ and $\gamma_5$ are nonparallel
in $S$, by \mbox{\fullref{cla:claim8}}, and none of them is parallel in $S$ to $\gamma_7$, by
\fullref{cla:claim9}. This implies that the genus of
$S$ is $ > 2$. 

So the third singularity must be of type 3. In this singularity the curves
$\gamma_1$ and $\gamma_3$ are fused into a curve $\gamma_5$; this is a trivial curve bounding in its 
interior two nonconcentric curves. Note that $\gamma_5$ must be essential in
$S$, for otherwise $S$ will be disconnected. If the next
singularity is of type 3 again, then $\gamma_2$ and $\gamma_4$ are fused 
into a curve $\gamma_6$, but then $\gamma_5$ and $\gamma_6$ must be parallel in $S$, contradicting 
\fullref{cla:claim8}.

If the fourth singularity is of type 4, then it can be interchanged with the third singularity,
showing then that the surface is compressible. If the fourth singularity is of type 5, the curve 
$\gamma_5$ is split into two curves which will be essential in $S$, but then there will be more than 3
essential nonparallel curves in $S$, which is not possible.

So the fourth singularity must be of type 2, splitting
$\gamma_2$, say, into $\gamma_6$ and $\gamma_7$.  The next singularity must be of type 4, fusing
$\gamma_6$ and $\gamma_7$ into a curve $\gamma_8$, for otherwise $S$ will be compressible. 
Note that the curves $\gamma_3$, $\gamma_4$, $\gamma_5$, $\gamma_6$ and $\gamma_7$ are essential in $S$,
and that $\gamma_3$ and $\gamma_4$ are nonparallel in $S$. If $\gamma_8$ is nontrivial in $S$, then
$S$ is disconnected or have genus $> 2$. Then $\gamma_8$ is trivial and $\gamma_6$, $\gamma_7$
are parallel in $S$. We must also have that $\gamma_4$ and $\gamma_5$ are parallel in $S$.
So we must have one more singularity of type 4, fusing $\gamma_4$ and $\gamma_5$ into a trivial
curve $\gamma_9$. Then we can assume that $\gamma_8$ and $\gamma_9$ lie on $T_0$, and bound
nested disks of $S_0$ enclosing just one of the arcs of $k_0$.
This shows that $S$ is a surface of type 6.

Summarizing, we have the following case:

(a)\qua The sequence of singularity types is 1, 1, 3, 2, 4, 4, $S_0$ consists of nested two trivial disks, and $S$ 
is a surface of type 6.

Suppose now that $S_1$ consists of two nested disks which enclose just one of the two arcs of $k_1$.
A similar argument shows that these are the possible cases for the sequences of types of singularities:

(b)\qua 2, 2, 4, 1, 3, 3, $S_0$ consists of two nested disks, and $S$ is a 
surface of type 6, 
but in an inverted position, ie changing the roles of $R_0$ and $R_1$.

(c)\qua 5, 5, 6, 7, 3, or 5, 5, 7, 6, 3, $S_0$ consists of a disk and an annulus, 
which are nested, and $S$ is a surface of type 7, but in an inverted position. 

(d)\qua 5, 6, 5, 6, 3, or 5, 5, 6, 6, 3, or 5, 6, 2, 4, 3, $S_0$ consists of a trivial disk and $S$ is a 
surface of type 4, which looks inverted.

(e)\qua 5, 5, 6, 7, 8, 8, or 5, 5, 7, 6, 8, 8, or, 5, 7, 5, 6, 8, 8, or 2, 5, 5, 6, 8, 8, or 5, 7, 2, 4, 8, 8, or 
2, 5, 2, 4, 8, 8, and $S_0$
consists of two meridian disks and $S$ is a surface of type 5, which looks inverted.
\end{proof}

\begin{cla}\label{cla:claim18} Suppose $S_1$ consists of just one annulus. Then $S$ is a surface of type 1 or 8.
\end{cla}

\begin{proof} The annulus in $S_1$ then determines an annulus in $T_1$ which contains both arcs of $k_1$. 
An argument as above shows that
there are several possibilities for the sequences of types of singularities, these are: 

(a)\qua 6, 1, 3, or 7, 6, 3, $S_0$ consists of a trivial disk, and $S$ is a surface of type 1.

(b)\qua 6, 5, $S_0$ consists of an annulus, and $S$ is a surface of type 1. 

(c)\qua 7, 7, 3, $S_0$ consists of a trivial disk and an annulus, which are nested, and $S$ is a surface of type 8, which looks inverted. 
\end{proof}

\begin{cla}\label{cla:claim19} Suppose $S_1$ consists of two non-nested annuli 
(and nonisotopic to nested annuli). Then $S$ is a  surface of type 3, 8 or 9.
\end{cla}

\begin{proof} Each of the annuli on $S_1$ determines an annulus in $T_1$
which contains an arc of $k_1$. The possible sequences of types of singularities are:

(a)\qua 6, 6, 3, $S_0$ consists of a trivial disk, and $S$ is a surface of type 3.

(b)\qua 7, 7, 3, $S_0$ consists of a trivial disk and two annuli, which are nested, and $S$ is a 
surface of type 8, which looks inverted.

(c)\qua 6, 7, 3, or 7, 6, 3, $S_0$ consists of a disk and an annulus, which are nested, and $S$ is a 
surface of type 9, which looks inverted.
\end{proof}
 
\begin{cla}\label{cla:claim20} Suppose $S_1$ consists of two nested annuli. Then $S$ is a surface of type 4, 5, 7 or 8.
\end{cla}

\begin{proof} One possibility is that there is an annulus in $T_1$ containing both arcs of $k_1$, but by 
tracking the singularities as above, we can see that this case is not possible. So the arcs of $k_1$
are on different components of $T_1-S_1$. The possible cases for the sequence of singularity types are: 

(a)\qua 6, 6, 3, $S_0$ consists of a trivial disk, and $S$ is a surface of type 4.

(b)\qua 6, 7, 8, 8, or 7, 6, 8, 8, $S_0$ consists of two meridian disks, and $S$ is a surface of type 5, 
but it looks inverted. 

(c)\qua 7, 6, 3, or 6, 7, 3, $S_0$ consists of a trivial disk
and an annulus, which are nested, and $S$ is a surface of type 7, but it looks inverted.

(d)\qua 7, 7, 3, $S_0$ consists of a trivial disk and two annuli, which are nested, and $S$ is a surface of type 8, which looks inverted.
\end{proof}

\begin{cla}\label{cla:claim21} Suppose $S_1$ consists of an annulus and a disk, which are non-nested 
(and cannot be isotoped to be nested). Then $S$ is a surface of type 3.
\end{cla}

\begin{proof} Here the annulus in $S_1$ determines an annulus in $T_1$ which contains an arc of $k_1$, 
and the trivial disk determines a disk in $T_1$ which contains the other arc of $k_1$. The possible 
sequences of types of singularities are:

(a)\qua 6, 5, 6, 3, and $S_0$ consists of a trivial disk.

(b)\qua 6, 2, 4, 3, and $S_0$ consists of a trivial disk.
\end{proof}

\begin{cla}\label{cla:claim22} Suppose $S_1$ consists of one annulus and one disk, which are nested. 
Then $S$ is a surface of type 4, 5, 7, 8 or 9.
\end{cla}

\begin{proof} If the disk in $S_1$ determines a disk in $T_1$ containing
just one of the arcs of $k_1$, then the possible sequences of types of singularities are:

(a)\qua 6, 2, 4, 3, or 6, 5, 6, 3, $S_0$ consists of a trivial disk and $S$ is a surface of type 4, 
but it looks inverted.

(b)\qua 7, 2, 4, 8, 8, $S_0$ consists of two meridian disks, and $S$ is a surface of type 5, which looks inverted.

If the disk in $S_1$ determines a disk in $T_1$ containing both arcs of $k_1$ then the possible sequences of 
types of singularities are:

(c)\qua 1, 8, 5, or 1, 5, 8, $S_0$ consists of two nested annuli, and $S$ is a surface of type 7.

(d)\qua 1, 8, 5, 6, 6, or 1, 5, 8, 6, 6, $S_0$ consists of two nested annuli, and $S$ is a surface of type 7.

(e)\qua 1, 8, 8, $S_0$ consists of an annulus, and $S$ is a surface of type 8.

(f)\qua 1, 8, 5 , or 1, 5, 8, $S_0$ consists of two non-nested annuli, and $S$ is a surface of type 9.
\end{proof}

\begin{cla}\label{cla:claim23} If $\vert S_1 \vert \geq 3$, then $S$ is a surface of type 8 or 9.
\end{cla}

\begin{proof} An argument as in \fullref{cla:claim17}, shows that it is not possible that $S_1$ contain 3 or
more trivial  disks, or two disks and some annuli. If $S_0$ contains a trivial disk and two or more annuli,
an argument as in previous claims shows that the only possibility is that the trivial disks bounds a disk in
$T_1$ containing both arcs of $k_1$, and that the disk and the annuli are nested. There are two possibilities:

(a)\qua $S_1$ consists of a disk and $n$ annuli, which are nested, the sequence of singularity types is 1, 8, 8,
$S_0$ consists of $n$ nested annuli, and $S$ is a surface of type 8.

(b)\qua $S_1$ consists of a disk and $n$ annuli, which are nested, the sequence of singularity types is 1, 8, 5 
(or 1, 5, 8), $S_0$ consists of $n+1$ annuli, two of them are innermost, the others are nested around the 
innermost ones, and $S$ is a surface of type 9.

Suppose now that $S_1$ consists only of annuli. By similar arguments as in previous claims, it can be shown
that there are two possibilities:

(c)\qua $S_1$ consists of $n$ annuli which are nested. The innermost one determines an annulus in $T_1$ which 
contains one arc of $k_1$. The other arc of $k_1$ is between the second and third annulus. The sequence of
singularity types is 7, 7, 3, $S_0$ consists of a trivial disk and $n$ nested annuli, and $S$ is a surface of type 8, which looks inverted.

(d)\qua $S_1$ consists of $n$ annuli, two of them are innermost, each determining an annulus in $T_0$ containing an arc
of $k_1$, the other annuli are nested around both of the innermost annuli. The sequence of singularity types is
6, 7, 3 (or 7, 6, 3), $S_0$ consists of a trivial disk and $n-1$ annuli, which are nested, and $S$ is a surface of type 9, which looks inverted.
\end{proof}

This completes the proof of \fullref{thm:maintheorem}.

The same arguments as in the previous claims can be applied when $S$ is a surface of genus 1,
and in fact the arguments are simpler. Again, the surface $S$ can be divided into pieces $S_1$, $S_0$, 
and $\tilde S$, and it can be shown that $S_1$ consists of just one disk or one annulus. So we have
the following result.

\begin{thm}\label{thm:tori} Let $K$ be a $(1,2)$--knot and $S$ a genus 1 meridionally 
incompressible
surface in the complement of $K$. Then $K$ and $S$ can be isotoped so that
$S$ look as one of the surfaces constructed in \fullref{sec:tori}.
\end{thm}

\section[Knots which are not (1,2)-knots]{Knots which are not $(1,2)$--knots}\label{sec:not1-2knots}

We recall the construction of \cite[\S~$6$]{E2}, which produces tunnel number one knots whose complement
contain a genus 2 closed meridionally incompressible surface which does not bound a handlebody in $S^3$.

Let $K$ be a satellite tunnel number one knot in $S^3$, and let 
$S$ be the closed incompressible surface of genus 1 contained in the complement of $K$;
then $S$ divides $S^3$ into two parts, denoted by $M_1$ and $M_2$, where, 
say, $K$ lies in $M_2$. In fact, it follows from \cite{MS} (or \cite{E1}) that
$M_1$ is the exterior of a torus knot, $M_2-\Int N(K)$
is homeomorphic to the exterior of a $2$--bridge link and that a fiber of $M_1$ is
glued to a meridian of $M_2-\Int N(K)$. Let
$\beta=\beta_1\cup \beta_2$ be an unknotting tunnel for $K$, where $\beta_1$ is
a simple closed curve, and $\beta_2$ is an arc joining $K$ and $\beta_1$.
The tunnel $\beta$ can be chosen so that $\beta_1$ is disjoint from $S$, and that $\beta_2$ intersects
$S$ transversely in one point, so $\beta_1$ lies in $M_1$. 
The surface $S$ then divides $\beta_2$ in two arcs, $\beta_2'$ and $\beta_2''$,
where $\beta_2'$ joins $K$ and $S$, and $\beta_2''$ joins $S$ and $\beta_1$.

Let $\gamma$ be a simple closed curve contained in $\partial N(K)$.
Assume that the arc $\beta_2'$ connects $S$ with a point in $\partial N(K)$,
so that such a point lies on $\gamma$.
Consider the manifold $M=M_2-\Int N(K)$. This is a compact,
irreducible $3$--manifold, whose boundary consists of two incompressible
surfaces, $S$ and $\partial N(K)$. The curve $\gamma$ lies on $\partial N(K)$.
We assume that $\gamma$ goes at least twice longitudinally around $N(K)$, and that
$\gamma$ is a nonreducing curve for
$M$, ie $M(\gamma)$, the manifold obtained by doing Dehn filling on $M$ with slope
$\gamma$, has incompressible boundary.

Note that if $K$ is not a cable knot, then any such curve $\gamma$ is a
nonreducing curve, for it is at distance $\geq 2$ from a meridian of $K$ \cite{W2}.
If $K$ is a cable knot, then there is a properly embedded annulus in $M$, with one of its 
boundary components lying in $S$ and the other one lying in 
$\partial N(K)$, which we denote by $\gamma_1$. It follows from \mbox{\cite[Theorem 2.4.3]{CGLS}}
that $M(\gamma)$ is $\partial$--reducible if and only if 
$\Delta(\gamma,\gamma_1)\leq 1$. 

Let $H=M_1 \cup N(\gamma)\cup N(\beta_2')$. We can assume that $H$ is made
of the union of the solid torus $N(\gamma)$ and the manifold $M_1$, 
which are joined by the $1$--handle $N(\beta_2')$. 
Let $W=S^3-\Int H$.
It follows from \cite[Theorem 6.3]{E2} that
$\Sigma =\partial W = \partial H$ is incompressible in $W$.

Let $K^*$ be a knot such that 
$K^*\subset H$, $K^*=K_1\cup K_2$, where $K_1$ is an arc contained
in $\partial N(K) \cap N(\gamma)$, and $K_2$ is an arc 
contained in $N(\beta_2')\cup M_1$,
which is obtained by sliding $\beta_1$ over $\beta_2$, ie $K_2$ is
an unknotting tunnel for $K$. In other words, $K^*$ is an iterate of $K$ and $\beta$, as
defined in \cite[\S~$6$]{E2}, and in particular $K^*$ is a tunnel number one knot.
Assume further that the wrapping number of $K_1$ in $N(K)$ is $\geq 2$,
that is, if we connect the endpoint of $K_1$ with an arc lying in $N(\beta_2')\cap N(\gamma)$,
we should get a knot whose wrapping number in $N(K)$ is $\geq 2$.  
It follows from \cite[Theorem 6.4]{E2} that $\Sigma$ is meridionally incompressible 
in $S^3-K^*$. 

Note that the wrapping number of $\gamma$ in $N(K)$ is $\geq 2$,
and that the wrapping number of $K_1$ in $N(\gamma)$ is also $\geq 2$, so
the wrapping number of $K_1$ in $N(K)$ is $\geq 4$. As $K$ is a $(1,1)$--knot,
then it seems that $K^*$ is a $(1,n)$--knot with $n\geq 4$.

Here, we show the following.

\begin{thm}\label{thm:not1-2knots} A knot $K^*$ constructed as above is not a $(1,2)$--knot, 
ie $b_1(K^*)\geq 3$.
\end{thm}

\begin{proof} Note that the surface $\Sigma$ does not bound a handlebody in $S^3$. In fact, to one 
side it bounds the manifold $W$, which has incompressible boundary, and to the other side
it bounds the manifold $H$, which is the disk sum of a solid torus (ie $N(\gamma$)) and the exterior of a 
torus knot (ie $M_1$).
This shows immediately that $K$ is not a $(1,1)$--knot, for any meridionally incompressible
surface in the complement of a $(1,1)$--knot bounds a handlebody in $S^3$ \cite{E4}.
It follows also that $\Sigma$ cannot be a surface of type 1, 2, 3, 4, 5, 8 or 9 for any such
surface bounds a handlebody in $S^3$. So, if we show that $\Sigma$ cannot be a surface of type 6 or 7,
then we will show that $K^*$ cannot be a $(1,2)$--knot.

By construction $K^* \subset H$, and there is a disk $D_W$ properly embedded in $H$, which intersects $K^*$ in two
points, separates $H$ and $\partial D_W$ is essential in $\partial H$. The disk $D_W$  is just the cocore of
$N(\beta_2')$, and we can assume that $D_W$ divides $K^*$ into the arcs $K_1$ and $K_2$ defined above. We claim
that if $D$ is another disk properly embedded in $H$, intersecting $K^*$ in two points, separating $H$, 
and with $\partial D$ essential in $\partial H$, then $D$
must be isotopic to $D_W$. To see that, look at the intersections between $D_W$ and $D$, and by doing an innermost
disk/outermost arc argument, remove all curves and arcs of intersection. To do that we use the following facts: (a)
$\Sigma$ is meridionally incompressible in
$H-K^*$;  (b) the arc obtained by sliding $\beta_1$ over $\beta_2''$ is not isotopic into the surface
$S$, for it is an unknotting tunnel for $M_1$; (c) the knot in $N(\gamma)$ ($M_1$) obtained from 
the arc $K_1$ ($K_2$), by joining the endpoints of $K_1$ ($K_2$) lying in $D_W$ and then pushing it
into $N(\gamma)$ ($M_1$), does not have local knots in $N(\gamma)$ ($M_1$). If $D$ and
$D_W$ are disjoint, then it is not difficult  to see that they are isotopic in $H$.

Let $S_6$ be a surface of type 6, and let $K_6$ be a $(1,2)$--knot in the complement of $S_6$,
so that the surface $S_6$ is meridionally incompressible in $S^3-K_6$. It follows from \fullref{sec:type6} that 
$S_6$ bounds a manifold $M_6$ which has incompressible boundary (this is a single manifold). 
Let $H_6= S^3-\Int M_6$. It follows from \fullref{sec:type6}  
that $H_6$ is the disk sum of the exterior, $S^3-\Int N(\gamma_6)$, of a certain knot
$\gamma_6$, and the solid torus $N(\alpha_1)$. We assume that $\gamma_6$ is a nontrivial knot, for otherwise
$H_6$ is a handlebody. Here there is also a disk $D_6$, separating $H_6$, with $\partial D_6$ essential in
$\partial H_6$,  and which intersects $K_6$ in two points. This is just a cocore of $N(\alpha_2)$. Note however 
that the disk $D_6$ may not be unique, it will depend on the way the corresponding arc of $K_6$ is embedded in
$S^3-\Int N(\gamma_6)$.

Let $S_7$ be a surface of type 7, and let $K_7$ be a $(1,2)$--knot in the complement of $S_7$,
so that the surface $S_7$ is meridionally incompressible in $S^3-K_7$. It follows from \fullref{sec:type7}
that $S_7$ bounds a manifold $M_7$ which has incompressible boundary (where this is a family of manifolds
constructed in a similar manner). Let $H_7= S^3-\Int M_7$. It follows from \fullref{sec:type7} that $H_7$ 
is the disk sum of the solid torus $N_1$ and the exterior  of a torus knot ($S^3-\Int N_2$), which we are
assuming is nontrivial. Here there is also a disk $D_7$, separating $H_7$, with $\partial D_7$ essential
in $\partial H_7$, and which intersects $K_7$ in two points. This is just a cocore of $N(\alpha)$. Note however 
that the disk $D_7$ may not be unique, it depends on the way the corresponding arc of $K_7$ is embedded in
$S^3-\Int N_2$.

If $K^*$ is a $(1,2)$--knot, there must be a homeomorphism $h: (S^3,K^*)\rightarrow (S^3,K_i)$, $i=6$ or $7$,
where $K_i$ is a $(1,2)$--knot having a surface of type 6 or 7. The image of the surface $\Sigma$ must be
a surface of type 6 or 7, and the image of $W$ must be a manifold $M_6$ or $M_7$. But in
the complement of $W$ there is a unique disk $D_W$ with certain properties, so in $M_6$ or $M_7$ there must be
also such a disk, ie the disk $D_6$  or $D_7$ must be the unique disk with that properties.
Then the homeomorphism must take the disk $D_W$ onto a disk parallel to $D_6$ or $D_7$.
This implies that $W\cup N(D_W)$ must be homeomorphic to $M_7\cup N(D_6)$ or to $M_7\cup N(D_7)$.

Note that $W\cup N(D_W)$ is just the manifold $M_2-\Int N(\gamma)$, which is an irreducible manifold
with incompressible boundary, and not homeomorphic to $T^2\times I$, $T^2$ a torus.
But $M_6\cup N(D_6)$ is homeomorphic to $N(\gamma_6)-\Int N(\alpha_1)$, which is a reducible manifold
with compressible boundary. And note that $M_7\cup N(D_7)$ is the manifold $(N_2-\Int (N_1\cup N(\alpha))\cup
N(\alpha)$, which is just $N_2-\Int N_1$, and it is homeomorphic to $T^2\times I$. 
So we got different manifolds in each case.

We conclude that $\Sigma$ cannot be a surface of type 6 or 7.
\end{proof}

\begin{figure}[ht!]
\centering
\includegraphics[scale=.75,angle=90]{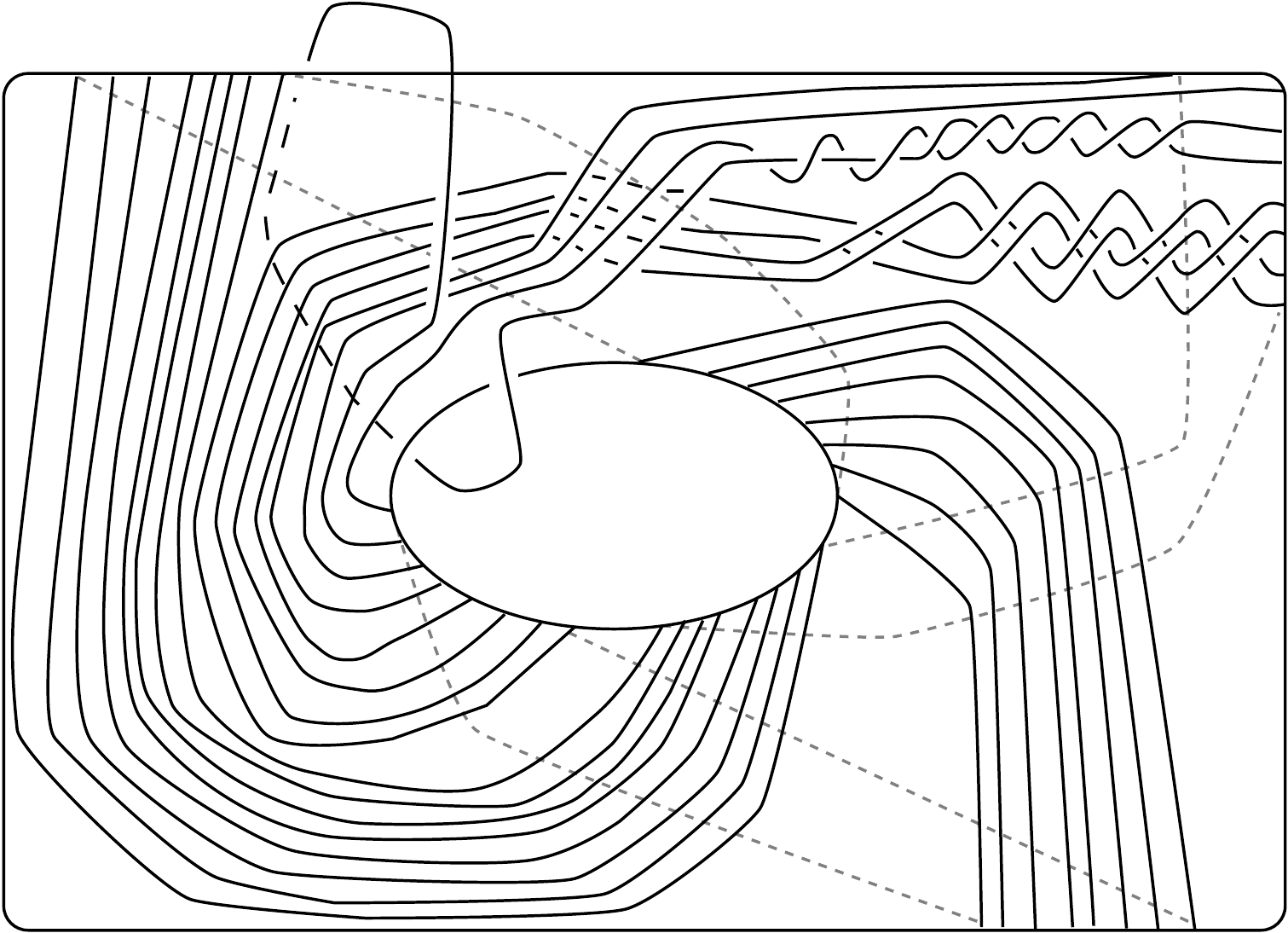}
\caption{}\label{fig:not1-2knot}
\end{figure}

\begin{figure}[ht!]
\labellist
\small\hair 2pt
\pinlabel $\beta_2'$ at 128 214
\pinlabel $\gamma$ at 69 226
\endlabellist
\centering
\includegraphics[scale=.7]{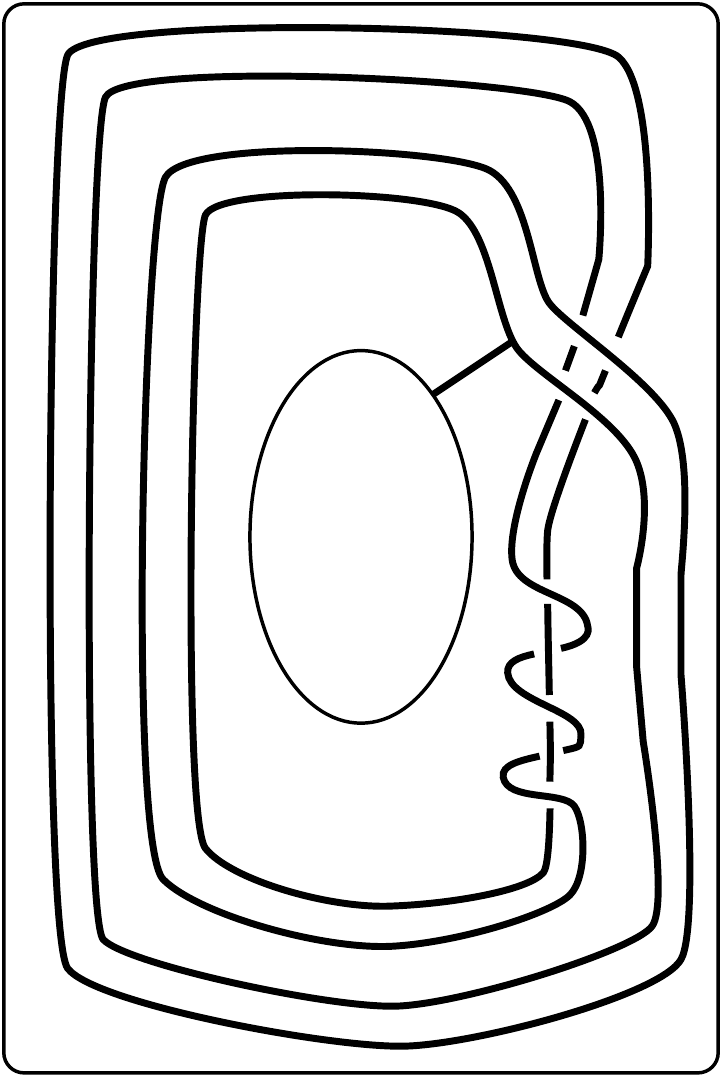}
\caption{}\label{fig:knotsurface}
\end{figure}

We give now an explicit example of a knot $K^*$. Suppose that $K$ is the $(-13,2)$--cable
of the left hand trefoil knot. This is a tunnel number one knot. The cabling annulus of $K$
has slope $(-26,1)$ on $K$. So let $\gamma$ be a curve on $\partial N(K)$ of slope $(-55,2)$.
Now take an arc $K_1$ on $\partial N(K)$ which goes around $N(\gamma)$ twice, and connect it with 
an unknotting arc $K_2$ for $K$. Such a $K^*=K_1\cup K_2$ is shown in 
\fullref{fig:not1-2knot}.
As said above, there is a surface $\Sigma$ in the complement of $K^*$ which is meridionally incompressible;
this surface is implicit in \fullref{fig:not1-2knot}. The surface $\Sigma$ bounds a manifold $W$, in this
example the manifold $W$ is homeomorphic to the manifold shown in 
\fullref{fig:knotsurface}, ie the exterior
of $\gamma\cup \beta_2'$ in the solid torus. To see a real picture,
just embed $W$ appropriately in the neighborhood of a trefoil knot, or of any torus knot.

\subsubsection*{Acknowledgments}

I want to thank Enrique Ram\'{\i}rez-Losada for several stimulating conversations.
I am grateful to the Technion in Haifa, Israel, for its support 
during the workshop held in July 2005, and specially to Yoav Moriah, who suggested me
several times to prove that there are tunnel number one knots which are not $(1,2)$--knots. 
I am grateful to the anonymous referee, who made a careful revision and
gave many helpful suggestions.
Research partially supported by PAPIIT--UNAM grant IN115105.

\bibliographystyle{gtart}
\bibliography{link}

\end{document}